\documentclass[10pt]{article}
\usepackage[top=2.0cm, bottom=1.5cm, inner=2.0cm, outer=2.0cm]{geometry}
\usepackage{color,graphicx,colortbl}
\usepackage{paralist,parskip}
\usepackage{booktabs, multirow, makecell}
\usepackage{subfigure, bm, comment}
\usepackage{amsmath,amsfonts,amssymb}
\usepackage{mathabx,fixmath,algorithm, algorithmic}
\usepackage{pstool}
\usepackage{psfrag}
\usepackage{hyperref}
\usepackage{cite}
\usepackage{natbib}

\definecolor{Gray}{gray}{0.9}

\DeclareMathOperator*{\minimize}{\textrm{minimize}}

\usepackage{theorem}
\newtheorem{remark}{Remark}

\newtheorem{definition}{Definition}

\input{mysymbol.sty}
\usepackage{needspace}

\newcommand{\rev}[1]{{{\color{black}#1}}}
\newcommand{\revtwo}[1]{{{\color{black}#1}}}
\newcommand{\revthree}[1]{{{\color{black}#1}}}

\title{\revtwo{Real-time City-scale Ridesharing via Linear Assignment Problems}}
\author{Andrea~Simonetto, Julien~Monteil, and Claudio Gambella \\ \normalsize{\emph{IBM Research Ireland, Dublin, Ireland}}\\ \small\emph{emails: }\texttt{andrea.simonetto@ibm.com}, \texttt{julien.monteil@ie.ibm.com}, \\  \small \texttt{claudio.gambella1@ie.ibm.com}}

\date{}
 
\begin{document}

\setlength{\parskip}{2mm} 
\setlength{\parindent}{0pt}

\maketitle

\begin{abstract}

In this paper, we propose a novel, computational efficient, dynamic ridesharing algorithm. The beneficial computational properties of the algorithm arise from casting the ridesharing problem as a linear assignment problem between fleet vehicles and customer trip requests within {a federated optimization architecture}. \revtwo{The resulting algorithm is up to four times faster than the state-of-the-art, even if it is implemented on a less dedicated hardware, and achieves similar service quality.} Current literature showcases the ability of state-of-the-art ridesharing algorithms to tackle very large fleets and customer requests in almost near real-time, but the benefits of ridesharing seem limited to centralized systems. Our algorithm suggests that this does not need to be the case. \revtwo{The algorithm that we propose is fully distributable among multiple ridesharing companies.} By leveraging \rev{two datasets, the New York city taxi dataset and the Melbourne Metropolitan Area dataset,} we show that with our algorithm, real-time ridesharing offers clear benefits with respect to more traditional taxi fleets \revtwo{in terms of level of service}, even if one considers partial adoption of the system. In fact, \rev{e.g.}, the quality of the solutions obtained in the state-of-the-art works that tackle the whole customer set of the New York city taxi dataset is achieved, even if one considers only a proportion of the fleet size and customer requests. This could make real-time urban-scale ridesharing very attractive to small enterprises and city authorities alike. \rev{However, \revtwo{in some cases, e.g., in multi-company scenarios where companies have predefined market shares, we show that the number of vehicles needed to achieve a comparable performance to the monopolistic setting increases, and this raises concerns on the possible negative effects of multi-company ridesharing.}}
\end{abstract}

\rev{
}

\textbf{Acknowledgement:} This project has received funding from the European Union's Horizon 2020 research and innovation programme under grant agreement No 731993 (Autopilot). 


\vskip2cm


\section{Introduction}\label{sec:intro}


It is estimated that by 2050, 66\% of the population will live in cities~\citep{UniNationsESA2014}. One key challenge stemming from urban-centric societies is the one of urban mobility. It is believed that on-demand ridesharing services are a way to offer tremendous benefits in terms of alleviating traffic congestion and pollution as they essentially provide a way to size down taxi fleets and car ownership. By introducing high flexibility on trips, times, and leveraging the shareability of the journey~\citep{Santi2014}, ridesharing services could provide timely and convenient transportation to virtually anybody at virtually anytime, using a minimal number of cars. Calculations performed with the New York City taxi data~\citep{NYCdata} have shown that 20\% of the current 4 seat taxi cabs are enough to serve 98\% of the current trip demand~\citep{Alonso-Mora2017}, and it was suggested in~\cite{Zhang2015} that an on-demand mobility service in the city of Singapore would meet the personal mobility needs of the \emph{entire population} with only 40\% of the current vehicle fleet. This could have disruptive consequences and change cities and mobility drastically. In the meantime, ridesharing systems are also becoming business cases and companies like Uber, Lyft, and Via are already investing large capitals. A recent research report warns that current ridesharing systems are contributing to the growth of vehicles miles traveled in cities~\citep{UCDavis2017}. In this context, it is critical to provide tools to the public authorities in charge of the policy decisions~\citep{CheAF17}, so that the opportunities offered by well-designed on-demand ridesharing systems are not missed. \revthree{Autonomous vehicles could be beneficial for shared mobility systems, in the near future~\citep{NAZARI2018456, levin2017congestion, chen2016operations, MA2017124, KRUEGER2016343, FARHAN2018310}.} \revthree{In~\citep{LOKHANDWALA201845}, an agent-based model evaluates the
	environmental and energy benefits of ride sharing using traditional and autonomous taxis. The framework copes with the driver preferences and the different availability of autonomous and traditional taxis. The results on a New York City case study show that autonomous driving in ridesharing can potentially decrease the fleet size by up to 59\% with a limited increase of waiting time and travel distance.}

Working within the frame of on-demand ridesharing services, in this paper we develop a novel, computational efficient, and dynamic, ridesharing algorithm that provides solutions of comparable quality to the ones of~\cite{Alonso-Mora2017}, on the same data set, while being less computational demanding. Our algorithm is based upon a suitably constructed linear assignment problem between new customer trip requests and available vehicles in the taxi fleets, and fits within a federated optimization architecture, which enhances flexibility.  

\subsection{Related work}

With the massive adoption of smart phones and social technology, ridesharing services are operating in a dynamic, wide, and articulated environment, and facing a variety of challenges, hence making the research and development of such services a very hot topic in the transport community. For instance, some recent research investigates the advantages of introducing meeting points in terms of cost savings and congestion mitigation~\citep{StiASG15}, the consideration of riders' satisfaction and privacy rights~\citep{AivGHK16}\revthree{, the integration of ridesharing in multi-modal systems~\citep{YAN2018, LIU2018}}, the offering of tailored pricing schemes \revthree{\citep{SAYARSHAD2018192}} (e.g. for regular travelers such as commuters~\citep{LiuL17, MaZ17}), the study of the changes in travel patterns induced by such ridesharing systems~\citep{DonWLZ18}\revthree{, and the consideration of ridesplitting as a binary classification problem~\citep{chen2017understanding}}.

In this work, more specifically, we focus on dynamic ridesharing systems, where customers and their schedule are revealed over time. Our motivation is to provide a real-time ridesharing solution for urban-scale fleet sizes and customer requests. \rev{On the basis of the data available for New York City~\citep{NYCdata} and for the Metropolitan Melbourne Area~\citep{Meldata}}, this requires to deal with problems involving several thousands of vehicles for millions of trips to be scheduled daily, that is having access to algorithms that can deal with thousands of vehicles in a matter of seconds. Due to the inherently difficult underlying optimization problem~\citep{Agatz2012}, \rev{and the challenges posed by the real-time implementation}, the literature mainly focuses on finding \rev{scheduling policies, having a low probability of constraint violation}. With this in mind, interesting surveys for the dynamic case are the work of~\cite{Berbeglia2010} which overviews online and dynamic instances for hundreds and even thousands of vehicles when using insertion heuristics, the methodological work of~\cite{Agatz2012}, the collection of research papers on the topic from 2006 to 2015~\citep{FHA2015}, and the work of~\cite{Ritzinger2016a} which more specifically addresses dynamic and stochastic routing. As for technical papers, they tend to focus on rolling horizon approaches by recomputing the solution of a static ridesharing instance over time. The work of~\cite{Ascheuer2000} focuses on online ridesharing problems where no time-windows are considered and proposes different methods for a warm-start that incorporates the most current schedule. Insertion heuristic approaches for the dynamic case are explored in~\cite{Varone2014}. Models and algorithms for the dynamic ridesharing problem are presented in~\cite{Colorni2011, Nourinejad2014}, while stochastic variants are studied in~\cite{Schilde2011}.

Most of the aforementioned literature offers solutions that can achieve scheduling policies for hundreds to thousands of vehicles in a ``reasonable'' time frame, by using the most promising heuristics. As we shall present now, some recent works have circumvented the curse of dimensionality and scaled the computational time further down to near real-time by proposing practical algorithms based either on the assignment problem or on geographical/time proximity. In~\cite{Agatz2011}, the authors divide the ridesharing problem in a two phase procedure. First, one-to-one rider-driver matching and no time windowing is assumed, and the problem is formulated as an assignment problem, where travel times are minimized. A rolling horizon approach is then used to handle new customers. With this methodology, the authors can deal with 29k trips in about a minute of computing time. The authors also present more complex variations of the basic problem, handling the possibility to book round-trips, and the choice of being a rider or a driver. In these cases, the optimization problem becomes more complex and they use Lagrangian relaxations as heuristics.

\rev{In~\cite{Hosni2014}, Lagrangian Decomposition approaches are developed for taxi ride-sharing services. Heuristic methods address dynamic variants of the problem for instances with up $50$ taxis and $200$ passengers.} 
In~\cite{NouR16}, a decentralized agent-based model for dynamic ridesharing is proposed: the ride-matching options determined by multiple virtual agents are fed into a single-shot first-price Vickrey auction. The model is tested on a Sioux Falls network, which is composed of $24$ nodes, $76$ transportation links and $1.96$ million potential users of the service. The decentralized approach proves to be computationally more convenient than the centralized one, without impacting the level of service. Taxis schedules have been investigated in~\cite{Ma2013} for Beijing and in \cite{Tian2013, Huang2014} for Shanghai. For both cases, the proposed approach modifies the current schedule by inserting the new queries based on proximity. In~\cite{Ma2013}, the method uses geographical clustering, matching with shortest path-like algorithms, and scheduling. Since only small instances of customer-taxi matching and scheduling are considered per time-step, the method can handle over 30k taxis over a 3 months time frame. In~\cite{Tian2013, Huang2014}, the approach is similar but the dynamic-matching is based on integer programming. This method handles more than 10k taxis per day and roughly 500k trips. In~\cite{davidtorche}, \rev{a cluster-based matching algorithm, based on trip spatial features, is run prior to running the assignment problem to match drivers and riders. Instances with hundreds of vehicles in the Melbourne Metropolitan Area can be solved in real-time in a rolling horizon framework}. In~\cite{MASOUD201760}, a peer-to-peer mechanism to exchange rides is proposed and this proves to increase the proportion of customers served.

In~\cite{Alonso-Mora2017} (and subsequently in ~\cite{Alonso-Mora2017a}), the authors propose an on-demand ridesharing algorithm that can match trips with vehicles by leveraging the idea of shareability networks~\citep{Santi2014}. In particular, given a set of new requests, first, the algorithm constructs a pairwise request-vehicle graph encapsulating the rides that could be shared among new requests; second, it computes a graph of feasible trips and the vehicles that can serve them; third, it solves an integer linear program to match vehicles to trips; finally, it rebalances the remaining idle vehicles. By reducing the scheduling problem to a on-demand trip-to-vehicle matching, the authors are able to handle a fleet of 3000 vehicles, more than 3 million rides over a week, with near real-time assignments in a matter of seconds. One of the most important feature of the algorithm, with respect to the state of the art, is its ability to tackle ridesharing with vehicles of (theoretically) arbitrary capacity. When applied to the New York City taxi data, the algorithm shows that 2000 vehicles with a capacity of 10, or 3000 vehicles with a capacity of 4 can serve up to 98\% of the demand, with a mean waiting time of 2.8 minutes and a mean trip detour of 3.5 minutes.  


\revtwo{Finally, we make use of the basic concept of federated optimization which has emerged as a suitable candidate to handle large-scale optimization problems, where data is distributed among different computing devices~\citep{Konecny2016}. Federated optimization differs from classical distributed optimization in the way communications and computations are conceived. Generally speaking, in distributed optimization, the devices perform only simple computations, while communication between devices is essential. Conversely, in federated optimization, the devices are key in the computations and they perform as little communications as possible. Federated optimization stems from machine learning, where data is acquired by multiple devices, which perform local operations and do not share the \rev{raw} data, but a subset of their local results\footnote{In this paper, we borrow the concept of federated optimization for its meaning of devices, in our case vehicles, owning their own data and solving local optimization problems whose solution is then integrated in a centralized decision center in one shot. The original paper~\citep{Konecny2016} pertains machine learning, but we feel that federated optimization can be seen in a wider context, as here.}. Federated optimization closely resembles the edge computing paradigm which aims at optimizing cloud computing systems within the emerging internet-of-things architecture~\citep{shi2016edge}. In edge computing, the data processing occurs near the source of the data whenever possible so as to reduce the communications bandwidth.}


\subsection{Contributions of this paper}

First, our approach is similar but differs from the one proposed in~\cite{Alonso-Mora2017}. By enforcing that all new requests processed at the same time cannot be combined, i.e., a single vehicle can be matched only to one new request, we are able to construct a less computational demanding algorithm based on a linear assignment problem. The limitation of having only one new request per vehicle is not as restrictive as it may seem: new requests are batched and come in at regular sampling intervals, typically 10-30 seconds, so customers can be matched to already occupied vehicles in the next batch of requests. Intuitively, the loss due to the impossibility of combining requests during the same sampling period can be overcome by increasing the sampling frequency. \rev{We will elaborate this point further in Section \ref{sec:pf}. Yet, this simplification offers significant benefits: qualitatively similar performance,  much less computational effort, a plethora of already available tools for the linear assignment problems, such as very efficient algorithms, distributed and parallel algorithms.}

Second, our approach is built \revtwo{within} a federated architecture. On the cloud, the linear assignment problem is run for each batch of requests, and on the edge, vehicles can compute the cost induced by the addition of a new customer to their schedule. \rev{As we will explain, this novel approach is not just another decomposition scheme, and offers clear advantages in terms of flexibility}.  A context-mapping module is introduced to appropriately select the vehicles that could accommodate new customers based on geographic information, so that the capacitated vehicle routing problem~\citep{ralphs2003capacitated} to be solved becomes computationally affordable. This module also has the advantage of considerably reducing the communications requirements.

\revtwo{Third, as our approach is based upon a sequence of linear assignment problems, it is fully distributable among multiple companies \revthree{in a privacy-aware fashion} while preserving optimality of the assignment. We use a result of the literature related to the distributed auction algorithm, see~\cite{4739098,naparstek2014fully}. This is a particularly interesting feature to have as we can situate ourselves in the realistic situation of users having access to ridesharing mobility via their smartphones. Some users may be registered to the smartphone app of one company, others to another. In this situation, one just requires a central agent, for instance a city authority, to get access to the difference between the two best assignment costs of the requests (e.g. difference of best \revthree{assignment costs}) at each iteration of the assignment process.} \revthree{In the context of the rise of Mobility as a Service (MaaS) platforms, such an algorithmic feature is very appealing, since it does not require to disclose full information on the vehicles of the ridesharing companies.}

In this paper, we achieve a better computational efficiency \revtwo{(up to $4$ time faster)} than in the work of~\cite{Alonso-Mora2017}. There, the dynamic ridesharing algorithm was implemented using C++ and was running on a dedicated a 24 core 2.5 GHz machine, while having almost near real-time performance when dealing with the complete request pool of Manhattan~\citep{NYCdata}. We report a Python (2.7) implementation of the algorithm, and despite working with a denser network graph
consisting of 17,446 nodes compared with the 4,092 nodes considered in~\cite{Alonso-Mora2017}, \revtwo{achieves a reduction of computational time w.r.t.~\cite{Alonso-Mora2017} of more than $4$ times, delivering similar service quality. Our implementation can run comfortably in (near) real-time on a 1 core 2.7 GHz laptop for the complete request pool, thereby making ridesharing accessible to less dedicated hardware.} Furthermore, we are able to demonstrate that ridesharing is not only useful in the case of monopolistic economies, e.g., the ones that have access to the whole customer requests and the whole taxi fleet, or only in the case of massive adoption \rev{of vehicles}. Considering $5\%, 10\%, 20\%$ of the customer requests, we show that a respective \rev{percentage} 
 of the fleet size achieves similar gains in terms of the service rate, waiting and detour times as when all the customer requests and fleet are considered. This makes the ridesharing business appealing for start-up, small-medium enterprises, and city authorities alike. \rev{As a matter of fact, the service rate obtained with our approach is better than in~\cite{Alonso-Mora2017}, which may be due to the selection of a different cost function and to a higher sampling frequency, as we will argue. \revtwo{Apart from the fact that the method proposed in~\cite{Alonso-Mora2017} could be made faster (e.g., by bounding the computational time, or by reducing the number of requests per vehicle), our work demonstrates that a single-request assignment can perform better than a multi-request assignment, when the optimization is run frequently enough.} We further corroborate the numerical evaluation on a different dataset from the Melbourne Metropolitan Area, Australia.} 
 

In addition, city policy makers must ideally keep in mind that the total number of vehicles utilized by such ridesharing services must be limited as a function of the demand, in order to achieve the traffic and environmental benefits foreseen. \rev{We explore this concern by considering a two-company scenario, where the customer shares are allocated upfront. We show that in this case the number of vehicles to obtain a similar service rate than a monopolistic setting increases, thereby exposing possible negative effects of multi-company ridesharing.} 

In brief, the contributions we offer in this paper can be summarized as follows: 
\begin{itemize}
\item we propose a novel solution to the city-scale ridesharing problem that builds upon a federated architecture so as to maximize computational efficiency;
\item our solution recasts the ridesharing problem into a succession of batch processes that combines a linear assignment algorithm, a context-mapping algorithm and a capacitated vehicle routing problem with pick-up, delivery and time-window;
\item \revtwo{the resulting algorithm can be distributed among multiple companies, which do not need to share any knowledge on their own proprietary system (location of the vehicles, routing service);}
\item the solution provides comparable results to the ones in the literature in terms of the quality of service achieved, i.e., the percentage of requests satisfied given some waiting time constraints;
\item \revtwo{single-request assignments to vehicles provide very high service rates in a dynamic context, and suggest that an optimal myopic search for multi-request assignment is not necessary;}
\item one conclusion is that ridesharing services do not need to be monopolistic, as small market shares can lead to equivalently good results in terms of the quality of service achieved;
\item \rev{we illustrate via a simple example that special care from the city authorities is needed in order to avoid negative effects in some multi-company scenarios.} 
\end{itemize}

\subsection{Organisation of the paper}

The remainder of the paper is organized as follows. Section~\ref{sec:pf} reports the formulation of the problem. The federated optimization architecture approach is described in Section~\ref{sec:fed}, where all the components of the proposed algorithm are described along with their theoretical properties. Section~\ref{sec:num} presents the numerical implementation and results for the New York City taxi dataset \rev{and the Melbourne Metropolitan Area dataset}. We conclude with a discussion and suggestions on future research directions in Section~\ref{sec:open}.


\section{Problem formulation}\label{sec:pf}

Given a set $\mathcal{M}$ of $m$ customer trip requests at time $t$ and a ridesharing fleet $\mathcal{C}=\{1, \dots, n\}$ of vehicles, each with capacity $C_i$, we are interested in determining ridesharing solutions in real-time. Our ridesharing service aim at providing an assignment of requests (customers) to available vehicles and their correspondent routes, according to some optimization criteria. Available vehicles are those that can pick up customers, while complying with the time constraints associated with the requests, and without exceeding their seat capacity.

We express the ridesharing service as an optimization engine, which is run at specific time periods $t_k$ ($k =0,1,2,\dots$). At each such time instant $t_k$, the service processes the requests submitted by customers in the time window $(t_{k-1}, t_k]$, and find optimal vehicle-costumer assignment and correspondent routes. A trip request is defined formally as follows. 
\begin{definition} \emph{(Trip request)} \label{def:trip}
A trip request $r$ is a 5-tuple, containing the origin coordinates $O\in\mathbb{R}^2$, the destination coordinates $D\in\mathbb{R}^2$, the time-windows constraints for pick-up and delivery: $t_O = (\underline{t_O}, \overline{t_O})\in\mathbb{R}_{+}^2$, $t_D = (\underline{t_D}, \overline{t_D})\in\mathbb{R}_{+}^2$, where $\underline{t_O}$ and $\underline{t_D}$ are the earliest time at origin and destination, respectively, and $\overline{t_O}$ and $\overline{t_D}$ are the latest time at origin and destination, respectively, and additional time constraints conveniently stored as a vector of real numbers $\Omega$. Thus, we can write $r = (O, D, t_O, t_D,\Omega) \in \mathbb{R}^{q}$, with $q \geq 4$.
\end{definition}
Definition~\ref{def:trip} is rather flexible and enables to impose further request-related constraints. In the remainder of the section, we describe: the possible time limitations associated with requests, the sampling frequency of requests, a specific choice in the problem formulation, the mathematical formulation, and finally the solution method.

\paragraph{Time Limitations.} The vector $\Omega$ of Definition~\ref{def:trip} expresses time limitations of the trip requests, such as:
\begin{itemize}
	\item maximum waiting time $\delta$: difference between $\overline{t_O}$ and the current time;
	\item maximum detour time $\Delta$: upper limit on the extra trip time allowed when modifying the current route to serve other customers;
	\item maximum journey length $\Gamma$.
\end{itemize}
In an on-demand ridesharing system, the customer waiting time indicates the delay in the algorithm reactivity, while the detour time measures the flexibility of the solution approach.  Ridesharing algorithms trade-off detour times with occupancies and reduced fleet size. The impact of time limitations on the ridesharing solutions will be investigated in Section~\ref{sec:num}.



\paragraph{Sampling Frequency.} We assume that the sampling period for new requests (i.e., $t_k-t_{k-1}$) is constant over the service time. In order to limit waiting times of customers, the ridesharing service should cope with short sampling periods (on-line mode). However, solutions of higher \rev{(but nonetheless myopic)} quality can typically be obtained in case the sampling period is longer, and hence the requests are processed in batches. In this paper, we find a trade-off between on-line and batch mode by considering sampling times of the order of $10$ seconds.

\paragraph{Problem Feature.} 

\rev{

\begin{figure}
\centering
\includegraphics[width = 0.7\textwidth]{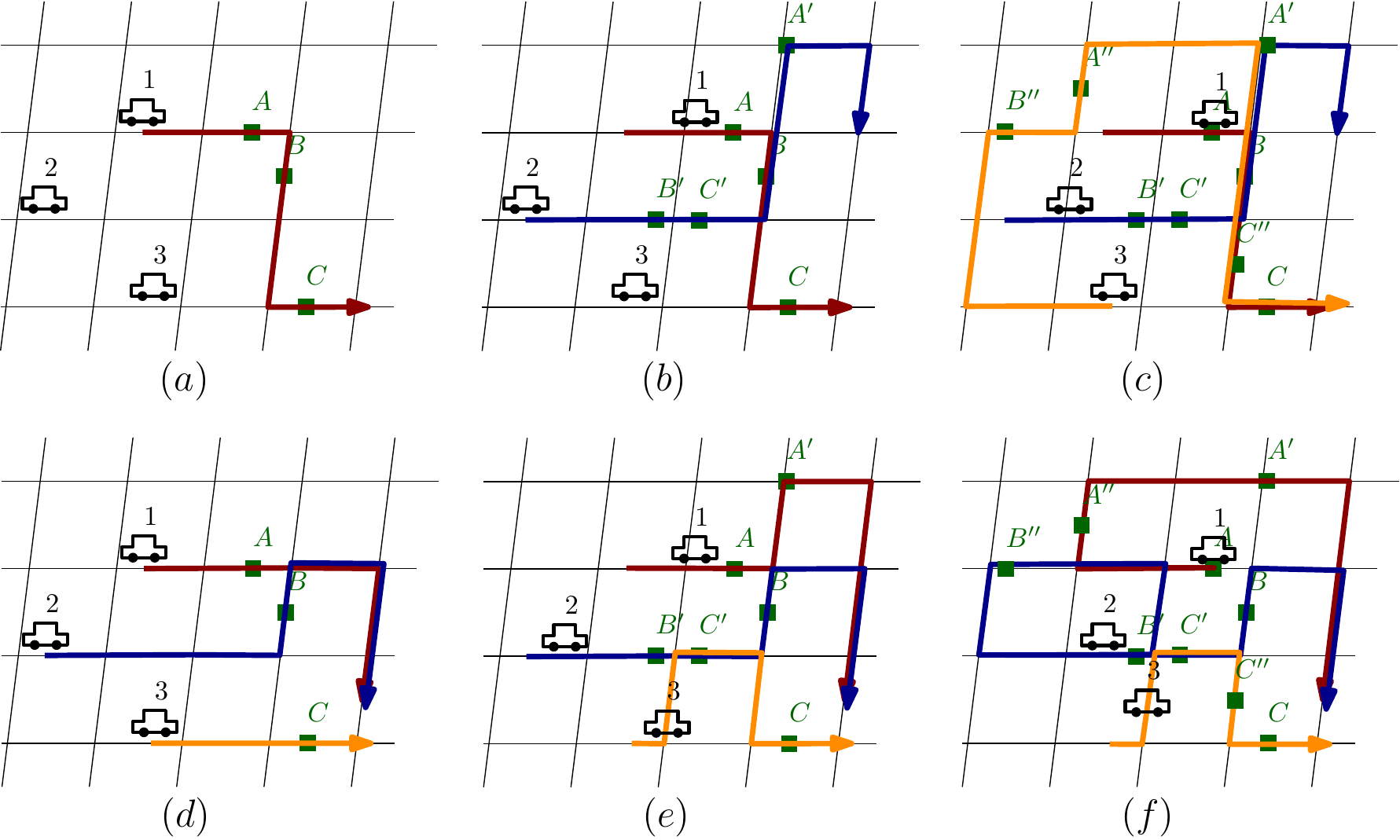}
\caption{Schematic of an hypothetical scenario of one-to-multiple assignment $(a)$-$(b)$-$(c)$ and one-to-one assignment $(d)$-$(e)$-$(f)$, for capacity $3$ vehicles.}
\label{fig:exp}
\end{figure}

Our algorithm leverages the following feature: at each optimization run, at most one new request (customer) is assigned to a vehicle. This design choice enables us to reduce the ridesharing optimization problem into a linear assignment problem, which can be solved efficiently and in a very fast manner. We consider this assumption reasonable in dynamic ridesharing problems. Let us consider the hypothetical scenario in Figure~\ref{fig:exp}: three vehicles with capacity $3$ are empty at time $0$. 

In $(a)$-$(b)$-$(c)$, we consider a one-to-multiple assignment, as the one in~\cite{Alonso-Mora2017}, where customers $(A,B,C)$, $(A',B',C')$, and $(A'',B'',C'')$ are revealed over three subsequent instances $t_1=10$~s $(a)$, $t_2=20$~s $(b)$, and $t_3=30$~s $(c)$. As we see, vehicle $1$ is associated with $(A,B,C)$ at $t_1$, vehicle $2$ with $(A',B',C')$ at $t_2$, and vehicle $3$ with $(A'',B'',C'')$ at $t_3$. 

In $(d)$-$(e)$-$(f)$, we consider a one-to-one assignment, as the one we consider here, where customers $(A,B,C)$, $(A',B',C')$, and $(A'',B'',C'')$ are revealed over three subsequent instances $t_1=10$~s $(d)$, $t_2=20$~s $(e)$, and $t_3=30$~s $(f)$. As we see, vehicle $1$ is associated with $(A,A',A'')$ over multiple times, vehicle $2$ with $(B,B',B'')$, and vehicle $3$ with $(C,C',C'')$.

What one can see is that:
\begin{itemize}
\item The one-to-one assignment choice does not prevent vehicles to serve multiple customers during the service;
\item Optimality at a given time does not guarantee the optimality at the end of the run: since both~\cite{Alonso-Mora2017} and our algorithm are myopic, it is not trivial to say which of the two strategy (one-to-multiple, one-to-one) is better in the long run. As one see from Figure~\ref{fig:exp}, assigning all the customers to one vehicle at the beginning, may prevent the same vehicle to serve a better suited customer that is coming later in time. 
\end{itemize} 

These considerations lead us to think that the one-to-one design choice is reasonable in dynamic ridesharing problems, and its benefits, such as the possibility to unlock fast and efficient algorithms for the assignment problem, over-weigh its defects. In practice, we see that our algorithm compares similarly to the one of~\cite{Alonso-Mora2017}, which copes with multiple customer assignments (and even performs re-assignment on-the-fly for not yet picked up customers).

}



\paragraph{Mathematical Formulation.} We introduce the mathematical formulation of the ridesharing optimization problem. Let $x_{ij}\in\{0,1\}, i\in\mathcal{C}, j \in \mathcal{M}$ be the set of binary variables: if $x_{ij} = 1$ then vehicle $i$ serves customer $j$, otherwise $x_{ij}=0$. When vehicle $i$ serves customer $j$, it incurs in a certain cost $c_{ij}\in\mathbb{R}_{+}$. The cost can be computed as a function of different parameters, typically: the detour time required to add the new customer to the current schedule, the vehicle capacity, the preference of the customers already onboard or scheduled to be serviced, etc. Such cost is denoted as a function $\mathcal{I}:\mathcal{C}\times\mathcal{M} \to \mathbb{R}_{+}$. In the case vehicle $i$ cannot accommodate the request of customer $j$, we set $\mathcal{I}(i,j) = \infty$, meaning infeasibility: this occurs when time windows or additional time limitations of $\Omega$ are not satisfied. The ridesharing optimization problem is then expressed as:
\begin{subequations}\label{eq:cas}
\begin{align}
\minimize_{x_{ij} \in \{0,1\}, i\in\mathcal{C}, j\in\mathcal{M}}\,& \hspace*{.5cm} \sum_{i=1}^n \sum_{j=1}^m c_{ij} x_{ij},  \label{obj}\\
 \textrm{subject to:}\, &\hspace*{.5cm} \sum_{i=1}^n x_{ij} = 1, \quad \, j\in\mathcal{M},\label{cons:requests}\\
 & \hspace*{.5cm} \sum_{j=1}^m x_{ij} \leq 1, \quad \, i\in\mathcal{C},\label{cons:vehicle}\\
& \hspace*{.5cm}c_{ij} = \mathcal{I}(i,j), \quad i\in\mathcal{C}, j\in\mathcal{M}.
\end{align} 
\end{subequations}
The aim of the optimization is to minimize the assignment costs \eqref{obj}: in our implementation, \rev{$c_{ij}$ is expressed as the time duration (or TD) of the route of vehicle $i$ that serves the already scheduled customers and the new customer $j$.}
Constraints \eqref{cons:requests} ensure that each request is assigned to a vehicle, while constraints \eqref{cons:vehicle} \rev{express that each vehicle can pick up at most a customer.} 
\revtwo{The practical implementation of the formulation is discussed in Section \ref{sec:optimization-module}.}

 \paragraph{Solution Method.} \rev{Provided that all the costs $c_{ij}$ are precomputed, then \eqref{eq:cas} can be casted as an asymmetric, or rectangular, linear assignment problem \citep{KenH80, PapS82, bijsterbosch2010solving}. It is well known that linear assignment problems have a totally unimodular constraint matrix, hence their continuous relaxation (i.e., substitute $x_{ij}\in\{0,1\}$ with $x_{ij}\in[0,1]$) is exact. We exploit the efficient solution algorithms for linear assignment problems (available for instances with up to $10^6$ requests \citep{Bernard2016, Bertsekas1989}). }
 





\section{Federated optimization architecture approach}\label{sec:fed}

\begin{figure}
	\psfrag{b1}{\scriptsize \hskip-.7cm Ridesharing logic}
	\psfrag{b2}{\scriptsize \hskip-.75cm Context mapping}
	\psfrag{b3}{\scriptsize \hskip-1cm Optimization module}
	\psfrag{b4}{\scriptsize \hskip-.2cm Vehicle }
	\psfrag{b1c}{\scriptsize \hskip-.2cm See Section 3.1}
	\psfrag{b2c}{\scriptsize \hskip-.2cm See Section 3.2}
	\psfrag{b3c}{\scriptsize \hskip-.2cm See Section 3.4}
	\psfrag{b4c}{\scriptsize \hskip-.2cm See Section 3.3}
	\psfrag{opt1}{\scriptsize \hskip-.5cm Assignment problem}
	\psfrag{opt2}{\scriptsize \hskip.5cm DARP}
	\psfrag{1-3}{\scriptsize \hskip-.2cm costs $c_{ij}$}
	\psfrag{3-1}{\scriptsize \hskip-.3cm assignment}
	\psfrag{3-o1}{\scriptsize \hskip-.3cm instance}
	\psfrag{o1-3}{\scriptsize \hskip-.3cm solution}
	\psfrag{1-2}{\scriptsize \hskip-.3cm trip request}
	\psfrag{2-1}{\scriptsize \hskip-.3cm vehicles}
	\psfrag{1-c}{\scriptsize \hskip-.3cm schedule }
	\psfrag{c-1}{\scriptsize \hskip-.3cm trip request }
	\psfrag{1-4}{\scriptsize \hskip-.3cm cost request }
	\psfrag{4-1}{\scriptsize \hskip-.3cm cost }
	\psfrag{o2-4}{\scriptsize \hskip-.3cm solution}
	\psfrag{4-o2}{\scriptsize \hskip-.3cm instance}
	\psfrag{Big-1}{Cloud}
	\psfrag{Big-2}{Edge}
	\psfrag{Routing}{\scriptsize \hskip-.45cm Routing service}
	\centering
	\includegraphics[width=0.65\textwidth]{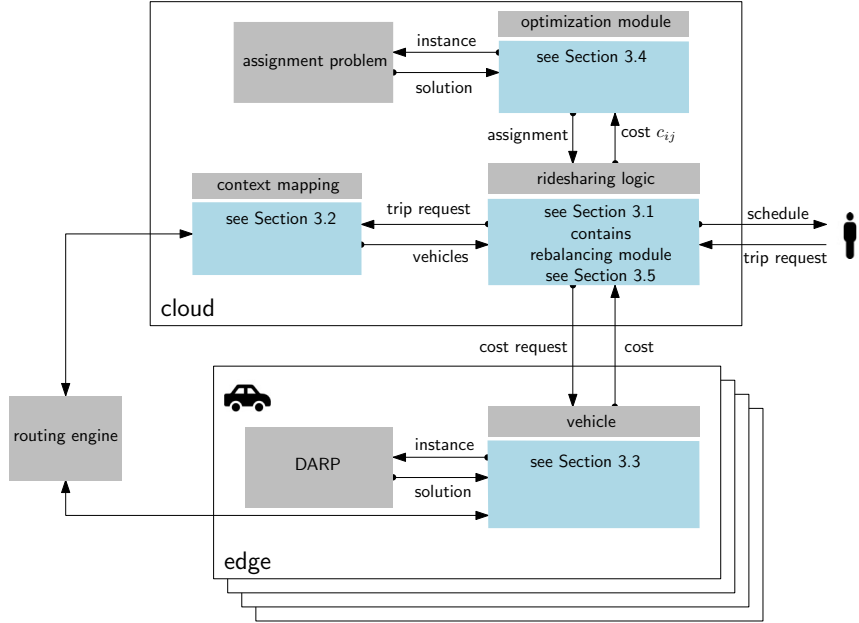}
	\caption{Overall dynamic ridesharing architecture. The cost $c_{ij}$ is linked to a specific schedule and route.}\label{fig.arch}
\end{figure}

Given the dynamic ridesharing problem at hand~\eqref{eq:cas}, the following federated structure is naturally identified:
\begin{itemize}
	\item \textbf{Bottom level}:
	\begin{itemize}
		\item \textit{Vehicle logic}: determination of the assignment costs $c_{ij}$ of customer $j$ to vehicle $i$. 
	\end{itemize}
	\item \textbf{Top level}:
	\begin{itemize}
		\item \textit{Context mapping}: algorithm for filtering the vehicles that could serve new customers, given their mutual geographic location. In other words, there is no reason to ask a vehicle on a route that is on \emph{that} side of the city to see whether it could accommodate a new customer on \emph{this} side of the city. This module considerably reduces the communications \rev{and computational} requirements.
		\item \textit{Optimization module}: given the assignment costs determined at the bottom level, this module provides a solution of Problem \eqref{eq:cas}.
		\item \textit{Ridesharing logic}: this block orchestrates the communication among the computational modules, as requests are processed.
	\end{itemize}
\end{itemize} 

The federated structure can be appreciated in Figure~\ref{fig.arch}. As in the federated optimization ``ideology'', the blocks solve local (almost private) optimization problems and share only limited pieces of information. The vehicle computational block can be replicated to accommodate as many vehicles as there are. The vehicle block needs access to a routing service (e.g., a webservice, or a dedicated GPS-navigation module): this is for determining the time-stamped route of each vehicle. The routing service is also queried by the context mapping, so as to compute the distances between vehicles and customers.

As one can infer from Figure~\ref{fig.arch}, the architecture structure is modular. Each vehicle could have their vehicle logic independently of each other, which means that their computation could be achieved in a parallel fashion. In addition, 
 both context mapping and assignment problem might be implemented on parallel architectures. While it is not so difficult to parallelize the context mapping module, for the assignment problem the interested reader is referred to~\cite{Wein1991}. Notwithstanding the parallel possibilities, in this paper (as mentioned in the introduction), we limit ourselves to single-threaded implementations. 

We describe the components of the federated optimization architecture approach in Sections \ref{sec:ridesharing-logic} - \ref{sec:optimization-module}. To limit the request rejection rate, we propose a reactive rebalancing procedure in Section \ref{sec:reactive-rebalancing}. Finally, the main theoretical properties of the ridesharing algorithm are discussed in \ref{sec:properties}.

\rev{Before moving on, a few considerations may be highlighted regarding the federated architecture structure. A federated optimization architecture is different from Lagrangian-based decomposition approaches applied to the overall optimization problem, even though the two approaches share similarities. In Lagrangian-based approaches, such as~\cite{Hosni2014,Mahmoudi2016}, the overall problem is automatically decomposed in simpler parts which are linked via Lagrangian multipliers to a master problem; such solution method provides a lower bound on the optimal solution value, because the assignments of requests to vehicles are not guaranteed to be feasible. The simpler parts perform local computations given a Lagrangian multiplier, while the master problem updates the multipliers. The process is iterative (back and forth from simpler parts to master problem) and eventually a solution is found. Lagrangian-based approaches have been demonstrated to be useful for limited fleet sizes and requests, e.g., in~\cite{Hosni2014} by using an incremental cost heuristic, the authors are able to tackle instances with 50 taxis and 200 passengers in 3.44 seconds, on average. 

The federated architecture presented in this paper can scale to much larger fleet size and requests. A federated architecture decomposes the problem less automatically, but tuned to the structure. The aim is to give as much computation as possible to the single physical entities (here the vehicles), which need to communicate as little as possible to a centralized server (ideally only once in both directions to find a solution). A federated architecture does not aim at finding the same solution as solving the overall problem would, but finding a feasible good solution with limited communication and time requirements. By communicating with the vehicles and back only once per request batch, the federated architecture reduces latency. By giving as much computation as possible to the vehicles, it enables them to make their own decision about the cost of adding a new customer to their schedule. Essentially, it means that any optimization problem and routing engine can be used, as this will not matter for the linear assignment solved at the centralized server. 
}

\subsection{Ridesharing logic}\label{sec:ridesharing-logic}

The ridesharing logic is run at every sampling period $t_k$ and works as follows:
\begin{enumerate}
	\item It obtains the customer requests submitted during the time interval $[t_{k-1}, t_k)$.
	\item It passes requests to the context mapping module. 
	\item It asks insertions costs to the vehicles returned by the context mapping module.
	\item It calls the optimization module for an optimal assignment.
	\item It sends to customers and vehicles the assignments and their corresponding routes. 
	\item \revtwo{If some customers cannot be serviced, it calls an internal rebalancing module, which runs the logic again from (2. to 5.) with loosen time constraints and for idle vehicles only. }
\end{enumerate} 
The ridesharing logic has access to a database containing the customers data (vector $r$ of Definition \ref{def:trip}) and vehicles data (e.g., current status, scheduled route). 


\subsection{Context mapping}\label{sec:context-mapping}


The first computational block is the context mapping block. This block is responsible to determine which vehicle should be asked for the computation of the cost $c_{ij}$ for each new customer $j$. \rev{

Given a new customer $j$, the context mapping:
\begin{enumerate}
\item Filters the available vehicles (of which are the ones that can accommodate a new customer in their schedules, which $n_i$ are idle and $n_o$ have already scheduled customers);
\item Computes the Euclidean distance between the available vehicles' position and the pick-up points;
\item Selects at most $2$max$n$ vehicles per request (max$n$ being a scalar):  the closest min$\{$max$n,n_i\}$ idle vehicles, and min$\{$max$n,n_o\}$ randomly picked already occupied (yet available) vehicles. This last step reduces the computational complexity. Note that, while it is clear that the closest idle vehicles are best in term of waiting times, the suitability of occupied vehicles  depends also on their schedule, so we pick them randomly. 
\end{enumerate}
}

Once the context mapping block has determined which vehicles should be asked to compute the cost of inserting which new customers in their schedule, the information is passed to the vehicle block. 

{\bf Computational requirements.} \revtwo{It is required to filter the vehicles, i.e., read their availability, which can scale linearly with the number of available vehicles and new customer requests, i.e., $O(n m)$; it is required to find the best $\max n$ elements (or pick $\max n$ elements at random)  in a list of at worst $n$ vehicles, with worst case complexity of $O(m n \log\, $max$n)$; finally it is required to compute Euclidean distances, which scale linearly with the number of available vehicles and new customer requests, i.e.,  in the worst case $O(n m)$.}

\subsection{Vehicle}\label{sec:vehicle}

The second computational block resides in the vehicle. It performs the determination of the vehicle-to-requests assignment costs. As discussed in Section \ref{sec:pf}, \rev{this requires to estimate time duration of the route of vehicle $i$ that serves the already scheduled customers and customer $j$.}


 Hence, the vehicle block is required to solve a single-vehicle Dial-a-Ride-Problem (DARP) \citep{Cordeau2003, LIU2015267, HAME201111}\revthree{, \citep{markovic2015optimizing}, \citep{HO2018395}}, with the aim of minimizing the route duration. In particular, the input data for the DARP \rev{solved by vehicle $i$} are: 
\emph{(i)} the request $r_j$,
 \emph{(ii)} \rev{the information relative to the scheduled route $R_i$, given by the pickup and delivery locations scheduled (either on the vehicle or in the pipeline) and the time limitations of the scheduled requests},
  \emph{(iii)} the current location $M_i$ of the vehicle and its capacity $C_i$,
   \emph{(iv)} the matrix \rev{$\tau$} of the travel times between pick-up/delivery locations and location $M_i$. 
The output of the DARP is expressed as 
\begin{equation}\label{darp}
(R_{ij}, D_{ij}) = \mathrm{DARP}(r_j, R_i, M_i, C_i, \tau),
\end{equation}
where $R_{ij}$ is an optimal route for vehicle $i$ serving request of customer $j$, starting from $M_i$, 
\rev{and $D_{ij}$ is the time duration of $R_{ij}$. In fact, in our formulation $c_{ij} = D_{ij}$.}

For vehicles (\rev{with no more that $3$ scheduled customers}), we solve Problem ~\eqref{darp} exactly by direct enumeration of the feasible insertion positions of the customer \rev{$j$}. \rev{This is because an exhaustive enumeration is tractable when the DARP deals with $4$ requests (i.e., the scheduled ones and the new one).} A pruning technique based on the current best solution limits the computational burden. For vehicles \rev{with capacity larger than $4$}, as in~\cite{Alonso-Mora2017}, we use \rev{the following insertion heuristic method (formalized in Algorithm~\ref{alg:ins-heur}): we start with the current schedule and route and attempt to insert the new customer pick-up and delivery, without changing the current schedule \emph{ordering}. So, if the current schedule is e.g., $(p1,p2,\ldots,d4,d1)$ meaning pick-up customer 1, pick-up customer 2, $\dots$, deliver customer 4, deliver customer 1, we try to insert p5 and d5 in all the possible positions that do not change the current ordering, i.e. $(p5,d5,p1,p2,\ldots,d4,d1), (p5,p1,d5,p2,\ldots,d4,d1), \ldots$. If an insertion is feasible, the total trip time is computed. At the end, the insertion that yields the lowest trip time is selected. }  

\begin{algorithm}[H]\caption{Insertion heuristic}
\begin{algorithmic}\label{alg:ins-heur}
	\STATE Set $\Lambda$ as the last position where to evaluate insertion of $r_j$
	\FOR{$k=1, \dots, \Lambda$}
		\STATE Let $\overline{k}$ be the $k$-th scheduled node in route $R_i$ not in location $M_i$.
		\STATE $R_{ij}^k := R_i$
		\STATE Add $r_j$ to $R_{ij}^k$ in position $\overline{k}$
		\IF{$R_{ij}^k$ is feasible for capacity and time windows constraints}
			\STATE $D_{ij}^k := $ time duration of $R_{ij}^k$
		\ENDIF
	\ENDFOR
	\STATE $k^* := \mathop{\mathrm{arg\,min}} \{D_{ij}^k\}_{k=1, \dots, \Lambda}$
	\STATE \textbf{return} $R_{ij} = R_{ij}^{k^*}, D_{ij} = D_{ij}^{k^*}$
\end{algorithmic}
\end{algorithm}

\rev{In Algorithm \ref{alg:ins-heur}, the possible insertion positions are limited by the parameter $\Lambda$, which is set to $4$ in our implementation. The insertion heuristic proves to yield good quality approximated solutions, as will be shown in the computational results. To evaluate the quality of Algorithm \ref{alg:ins-heur}, we implemented a Local Neighborhood Search (LNS) algorithm \citep{shaw1998using, pisinger2010large, HAME201111}) for vehicle fleets of limited-size.}
\rev{The solution of Algorithm \ref{alg:ins-heur} is fed into LNS, as described in Algorithm \ref{alg:lns}. Given the destroy and repair operators applied on vehicle routes, Algorithm \ref{alg:lns} is able to explore more possible insertions than those of the Algorithm \ref{alg:ins-heur}, and hence to improve the solution quality. The acceptance criterion for non-improving solutions is based on Simulated Annealing \citep{kirkpatrick1983optimization}.}
\begin{algorithm}[H]\caption{Local Neighborhood Search (LNS) for single-vehicle DARP}
	\begin{algorithmic}\label{alg:lns}
	\STATE Set $q \in \mathbb{N}$, $T_{start}$, $0<c<1$, $M_{iter}$
	\STATE Let $x$ be the solution obtained by Algorithm \ref{alg:ins-heur}
	\STATE Let $f$ be the function to evaluate the travel time of a route
	\STATE $x^{*}:= x$
	\STATE $T:= T_{start}$
	\WHILE{$i<M_{iter}$}
		\STATE $x_i := x$
		\STATE remove randomly $q$ requests from $x_i$ (\textit{destroy} $x_i$)
		\STATE reinsert removed requests into $x_i$, in a greedy fashion (\textit{repair} $x_i$)
		\IF{$f(x_i)<f(x^{*})$}
					\STATE $x^{*} \gets x_i$
				\ELSIF{$\textrm{accept}(x_{i},x, T_{start})$}
					\STATE $x \gets x_i$	
				\ENDIF			
		\STATE $T \gets T \cdot c$ 
		\STATE $i \gets i + 1$
	\ENDWHILE
	\STATE \textbf{return} $x^{*}$, $f(x^{*})$
	\end{algorithmic}
\end{algorithm}

\rev{Future research directions could investigate on more advanced heuristics (e.g., tabu search \cite{cordeau2003tabu} and adaptive large neighborhood search \cite{ropke2006adaptive}).

In our ridesharing problem, each vehicle is not a priori associated with a destination point, or a bound on route duration. Rather, it picks up and deliveries customers as long as necessary and compatible with seat availability. In order to simulate a limitation on the vehicle routes, we assume that a vehicle with capacity $C$ can have at most $2C$ requests in its pipeline of scheduled pickup/delivery events. This helps to lower the number of vehicle blocks computation at each time period.}


{\bf Computational requirements.} We compute the travel times \rev{$\tau$} via a routing server and store them offline; in this way, the computation burden of the vehicle block is the DARP solution. Said $N$ the number of customers and requests in the DARP instance, a complete enumeration of the feasible DARP solution has a computational complexity of  $O(\sqrt{N} (N^2/2)^N)$ where $N$ is the number of scheduled customers, if time-windows feasibility is neglected. Dynamic programming approaches \citep{Psaraftis1980}, are still non polynomial, having a complexity of $O(N^2 3^N)$. Our insertion algorithm has instead to evaluate only $O(N^2)$ candidate solutions. \rev{The LNS evaluates at most $N^2  M_{iter}$ solutions.}

\rev{
\begin{remark}
It is harder to say how the computational requirement differs in fleets of different capacity, because the real difference is the number of scheduled customers per vehicle. Consider two fleets of capacity $C=4$ and $C=10$ and suppose the number of scheduled customers per vehicle $C=4$ is a constant $\bar{N}$: the total cost would be $n O(\bar{N}^2)$. For $C=10$, one would have vehicles with more scheduled customers $\bar{N}_{+}$ and vehicle with less $\bar{N}_{-}$, say $n_1 O(\bar{N}^2_{+}) + n_2 O(\bar{N}^2_{-})$ with $n_1 +n_2 = n$. In general it is hard to say which complexity is the lowest.   
\end{remark}
}

\subsection{Optimization module}\label{sec:optimization-module}

\rev{The optimization module is responsible for the assignment of the customers to the vehicle routes. It proceeds as follows:
	\begin{itemize}
	\item it collects the assignment costs computed by the vehicle block;
	\item \revtwo{it adds virtual customers $\mathcal{V}_{\mathcal{M}}$ or artificial vehicles $\mathcal{V}_{\mathcal{C}}$ such that $\mathcal{M}'=\mathcal{M} \cup \mathcal{V}_\mathcal{M}$ and $\mathcal{C}'=\mathcal{C} \cup \mathcal{V}_\mathcal{C}$ have the same cardinality;}
	\item \revtwo{it sets to $\infty$ the assignment costs $c_{i, j}$ of artificial vehicles $\mathcal{V}_{\mathcal{C}}$ to all customers in $\mathcal{M}'$ and of real vehicles $\mathcal{C}$ to  virtual customers $\mathcal{V}_{\mathcal{M}}$}; 
	\item it equivalently reformulates the ridesharing problem \eqref{eq:cas} as a symmetric linear assignment problem \eqref{eq:lap}:
	\end{itemize}
	
}

\begin{subequations}\label{eq:lap}
\begin{align}
\minimize_{x_{ij} \in [0,1], \ i\in\mathcal{C}, j\in\mathcal{M}'}\,& \hspace*{.5cm} \sum_{i=1}^n \sum_{j=1}^n c_{ij} x_{ij}, \\
 \textrm{subject to:}\, &\hspace*{.5cm} \sum_{i=1}^n x_{ij} = 1, \, j\in\mathcal{M}',\\ 
&\hspace*{.5cm}  \sum_{j=1}^n x_{ij} = 1,\, i\in\mathcal{C}'.
\end{align} 
\end{subequations}
\rev{Considering a symmetric linear assignment problem enables to use highly efficient solution approaches, such as auction algorithms \citep{Bertsekas1981}.}
The solution of Problem \eqref{eq:lap} is the optimal assignment between the vehicles and the current batch of requests. The assignments \rev{contain} the specification of the routes (provided by \rev{the vehicle block and} the routing service), and those are sent to the ridesharing logic. \rev{In other words, the assignment solution prescribes: (i) the assignment of the satisfied requests to the vehicles, (ii) the position of insertion of the new customer in the vehicle route (computed in the DARP) and (iii) and the modified route (provided by the routing service). That is, one choses a particular assignment $ij$ based on the cost $c_{ij}$, which is linked to a specific insertion and route.  
	\revtwo{If an infinite cost is in the assignment solution}, then the customer cannot be accommodated and it is stored in a set $\mathcal{R}$ for further consideration in the rebalancing phase.} 

As of the solution approach, one can call any linear optimization solver available. Given the context mapping module, each customer might be assigned to a few vehicles in a time period, hence the assignment matrix of Problem \eqref{eq:lap} is sparse. We exploit the sparsity with a dedicated implementation of the auction algorithm~\citep{Bernard2016}. 


{\bf Computational requirements.} The computational complexity of state-of-art algorithms for linear assignment ranges between $O(n^3)$ and $O(n^3 \log(n))$. For the auction algorithms, the practical performance may scale much better~\citep{Bertsekas1981}.

\subsection{Reactive rebalancing}\label{sec:reactive-rebalancing}

To complete the ridesharing service, a rebalancing strategy is implemented (in a rebalancing module). This allows the fleet of vehicles to move towards highly demanded areas, while the algorithm is running. Similarly to~\cite{Alonso-Mora2017}, we implement a reactive rebalancing protocol. Whenever a new request (customer) cannot be accommodated, either because no vehicles is near them or because their time constraints are too tight, we consider all the idle vehicles and reroute one of them towards the customer. 

In practice, given a set of customers $\mathcal{R}$ that cannot be accommodated and would be refused and the set of idle vehicles $\mathcal{I}$, we determine which vehicles will be rerouted to which customer by solving the linear assignment problem
\begin{subequations}\label{eq:lap-r}
\begin{align}
\minimize_{x_{ij} \in [0,1], \ i\in\mathcal{I}, j\in\mathcal{R}}\,& \hspace*{.5cm} \sum_{i \in \mathcal{I}} \sum_{j\in \mathcal{R}} \rev{w}_{ij} x_{ij}, \\
 \textrm{subject to:}\, &\hspace*{.5cm} \sum_{i\in \mathcal{I}} x_{ij} = 1, \, j\in\mathcal{R},\\ 
&\hspace*{.5cm}  \sum_{j\in \mathcal{R}} x_{ij} \leq 1,\, i\in\mathcal{I},
\end{align} 
\end{subequations}
where $\rev{w}_{ij}$ is the trip time for vehicle $i$ to reach customer $j$ \rev{(calculated by a short-path algorithm and linked to a specific route)}. Problem~\eqref{eq:lap-r}, properly augmented to include virtual customers \revtwo{and artificial vehicles}, is a linear program as Problem~\eqref{eq:lap} and \rev{can be} easily solved. \rev{(Note that we have again a one-to-one assignment, while in~\cite{Alonso-Mora2017} one-to-multiple assignments, even in the rebalancing phase, are possible).}

It has to be noted that, depending on the number of available idle vehicles, \rev{some refused customers may fail to be accommodated by a rebalancing vehicle, and therefore they are refused}.

In addition, in our implementation, since the reroute vehicles do not typically satisfy the time constraints of the requests, the customers have the option to either accept the assignment (and perhaps be charged less), or refuse the match and be considered lost. 

\revtwo{Note that the rebalancing module is a part of the ridesharing logic and it calls once again the idle vehicles, for a different DARP problem with loosen time constraints, and an assignment problem. }

\subsection{Properties of the algorithm}\label{sec:properties}

Before assessing the algorithm on real data, we focus here on discussing its theoretical properties, in particular its computational complexity in relation with~\cite{Alonso-Mora2017}.

Computational complexity favorably compares with~\cite{Alonso-Mora2017}. The method of~\cite{Alonso-Mora2017} consists of four steps: (1) computing a pairwise request-vehicle shareability graph, (2)  computing a graph of feasible trips and the vehicles that can serve them, (3) solving an integer linear programming to compute the best assignment of vehicles to trips, (4) and rebalancing. In our algorithm, we perform similar steps, however: we do not compute the shareability graph (since requests cannot be combined in our formulation); step (2) has the same complexity of the vehicle computational block in our algorithm by computing which vehicle can be associated with which new request; step (3) is for us a linear assignment problem with significant lower complexity than the integer linear programming that is considered in~\cite{Alonso-Mora2017}; step (4) has a similar complexity of our rebalancing strategy (they are both linear programs with the same number of variables). In brief, our algorithm is computationally lighter than the one of~\cite{Alonso-Mora2017}.

The computational complexity of our ridesharing algorithm mainly lies in the computation of assignment costs and in the linear assignment problem.
As previously mentioned, the complexity of the cost computation scales linearly in the number of vehicles and the number of new requests, while worst-case complexity of DARP is non-polynomial in general, but polynomial for the insertion heuristic we use. As for the linear assignment problem, the number of variables in the worst-case are $n^2$, which happens only when all the vehicles are considered for the assignment problem. The complexity of the linear assignment depends on the algorithm used, being in the worst-case $O(n^3) - O(n^3 \log(n))$ depending on the algorithm. Practical performance (especially for auction algorithms) scales much better and, as mentioned, with dedicated implementations of the auction algorithm~\citep{Bertsekas1981}, solving sparse problems of the size of $10^6$ or dense ones of the size of $10^3$ as a matter of a few minutes. In our implementation, for a vehicle fleet of $3000$, the assignment problem takes approximatively from 3 to less than 5 seconds, due to the sparsity of the cost matrix. In comparison, the method for step (3) in~\cite{Alonso-Mora2017} uses heuristics to solve an integer non-convex optimization problem, with worst-case in the order of $O(n m^C)$ number of variables, whose scalability can be an issue.

\subsection{Distributed implementation}\label{distri}

\revtwo{Before moving on to the numerical assessment, we discuss the possibility to implement the algorithm in a distributed fashion, where multiple companies send limited information to a centralized server and the assignment problem is still solved to optimality. We note that to our knowledge, this possibility we offer is novel in the context of real-time city-scale ridesharing and is tied to the solution we have devised (since it hinges on distributed solutions of the linear assignment problem). \revthree{We also note that this approach appears very well suited to be implemented in MaaS platforms, as no proprietary information of the companies needs to be shared with the platform, and the limited information sharing still allows the linear assignment problems to be solved to optimality.}

Indeed, in a distributed scenario with multiple companies users will try to book a car via the proprietary applications of the companies. The context-mapping and cost-computation will then happen at the vehicle and company level. Then, as per the distributed auction algorithm (see Table 1 and Theorem 1 of~\cite{naparstek2014fully})   there is only a minimum amount of information that needs to be exchanged with a central agent/coordinator to setup and solve the assignment problem. 
The central agent needs to have access to the location of all customer trip requests, and their bids. In this way, the central agent cannot infer the location of all the vehicles of the companies. In detail, at each iteration the yet unassigned requests are computing their bids as the difference between the two lowest assignment costs (e.g., travel time differences) and share the bids with the central agent. Customers are then assigned to vehicles with the best bids and the assignment is communicated back to the companies. The process is repeated until all customers are assigned, and the optimal assignment is reached. Note that this allows each company to have their own proprietary system and not share any details about it, which is a much desired requirement in a realistic setting. The only requirement is a top-layer infrastructure, i.e., enforced by a city authority, that coordinates the bidding process.

The present paper implements the centralized auction algorithm. We will explore further the properties of this distributed architecture, \revthree{and its variants, among others,} in a future work. 
}

\section{Numerical implementation and results}\label{sec:num}

\subsection{Implementation}

We have implemented our ridesharing algorithm in Python 2.7 language on a 2.7GHz Intel i5 laptop with 8GB RAM memory. By doing this, we want to show that city-scale ridesharing solutions can be computed with non-dedicated hardware and software.


\rev{The algorithm was implemented using object-oriented programming. The ride-sharing logic, or federated scheduler, reads the stream of requests corresponding to one batch period and calls the context-mapping module that matches vehicles and requests. The context-mapping class takes the routing engine, the customer and vehicle data access objects to access the list of available vehicles and list of requests. Then the cost computation class is instantiated to compute the insertion costs of customers to vehicles making use of the DARP solver. The assignment problem is then solved via instantiation of the optimizer class which takes the cost matrix as input. Based on the solutions of the assignment problem, the status of the customers and routes of the vehicles are being updated, and the next batch of requests can be processed in turn.}

More specifically, the cost computation module implements a solver for the DARP problem in~\eqref{darp}, exact up to \rev{$3$ scheduled customers} and the insertion heuristic \rev{described in Algorithm \ref{alg:ins-heur}} for larger values of scheduled customers, while the time matrices \rev{$\tau$ are either extracted from a pre-computed node to node time matrix and stored using the Python package h5py (for the case of Manhattan), or computed on-line by the routing engine (for the larger case of the Melbourne Metropolitan Area)}. The routing engine relies upon calling a local OSRM (Open Street Routing Machine, \cite{Luxen2011}) server, which makes the generation of routes fast and efficient. \rev{For vehicle fleets with $100$ and $150$ vehicles with capacity $10$, we tested the LNS described in Algorithm \ref{alg:lns} with $M_{it}=10$ iterations, destruction degree $q=4$, initial temperature $T_{start}=10$ and cooling factor $c=0.9$.} For the linear assignment problem \eqref{eq:lap} we use the Matlab implementation of~\cite{Bernard2016} called via the Python Matlab engine, which scales better for this particular problem than general-purpose optimization solvers such as CPLEX \citep{Cpl} (when called via Python).

\subsection{New York City simulations}

We consider the New York City taxi public dataset~\citep{NYCdata} and extract one week of data trips from 00:00 hours on Sunday, May 5, 2013 to 23:59 hours on Saturday May 11, 2013 (exactly as in~\cite{Alonso-Mora2017}). This dataset contains the time and location of all of the pick-ups and drop-off locations visited by the 13,586 active taxis, for each day. From these data, we extract all of the origin-destination requests, together with the time of request equal to the time of pickup. Each day contains from 382,779 (Sunday) to 460,700 (Friday) requests. We consider the complete road network of Manhattan as encoded in Open Street Map, amounting at 17,446 nodes versus the 4,092 considered in~\cite{Alonso-Mora2017}.

The fleet is initialized in random locations within Manhattan at 22:00 hours on Saturday, May 4, 2013, (and the requests between 22:00 and 23:59 are used to warm start the fleet, so that the vehicles do not start empty at 00:00 on May 5) while the requests are batched and sent to the ridesharing service every $10$ seconds. We consider that each trip has constraints on the maximum delay $\delta$ and detour time $\Delta$, with $\delta = 7$ minutes, and $\Delta = 7$ minutes.

\begin{figure}
\centering
\includegraphics[width=0.55\textwidth]{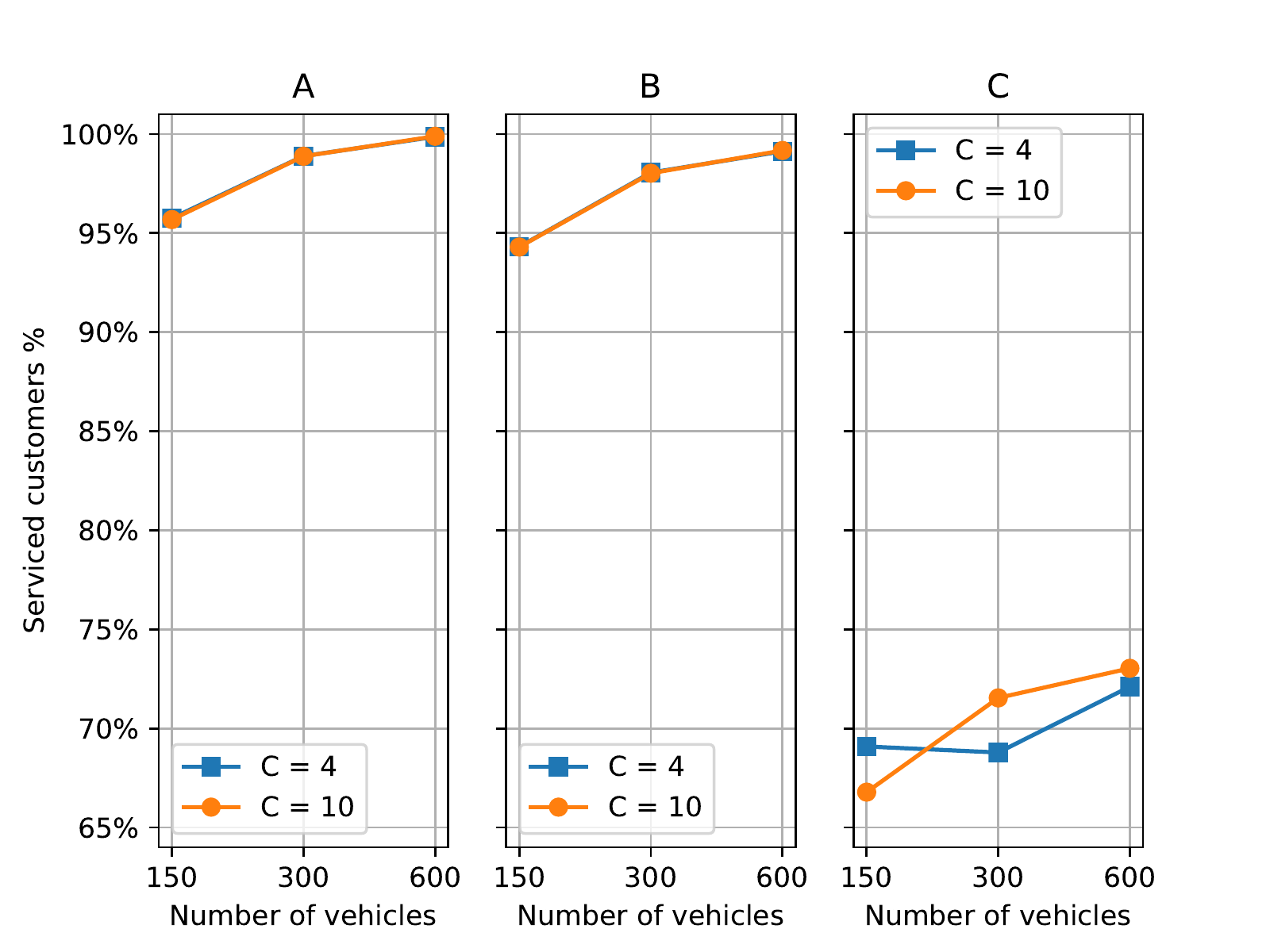}
\caption{Serviced customers: (A) with customers accepting the rebalancing vehicles, (B) with customers not accepting the rebalancing vehicles, (C) without rebalancing.}
\label{fig.sim-1}
\end{figure}

\begin{figure}
\centering
\includegraphics[width=0.99\textwidth]{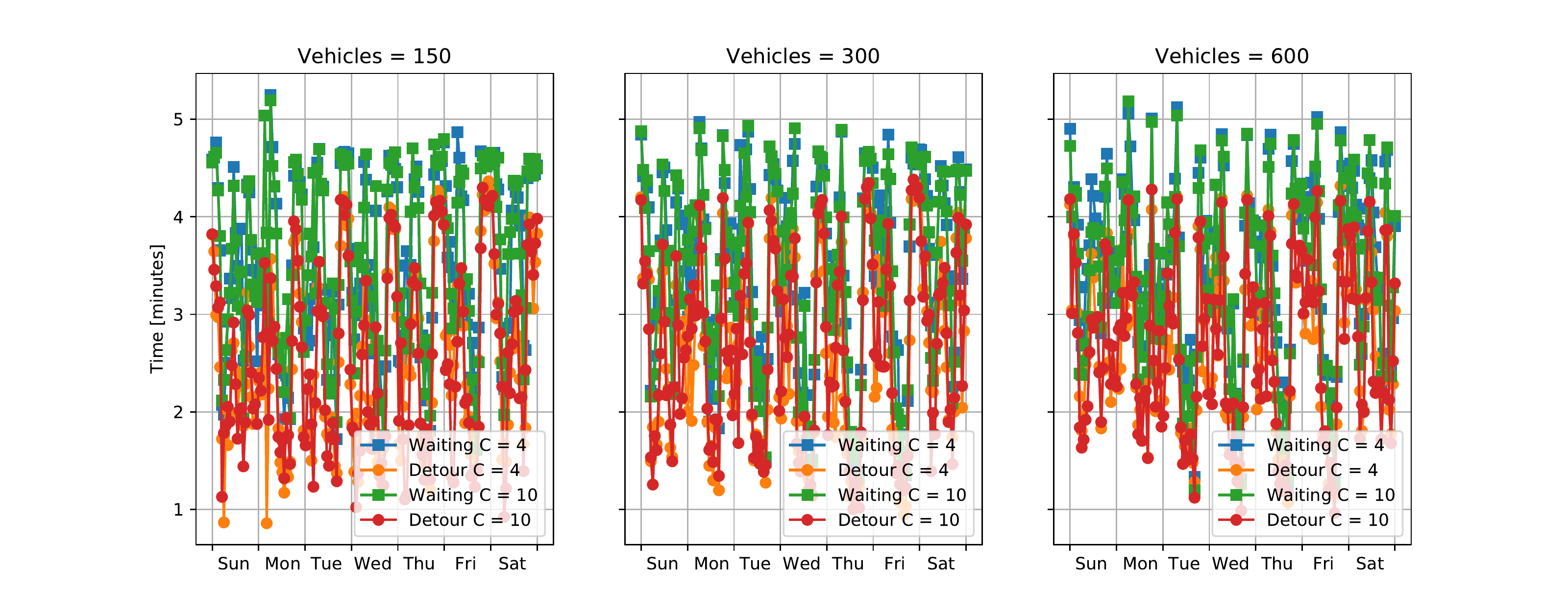}
\caption{Waiting and detour times per customers in the considered time frame.}
\label{fig.sim-2}
\end{figure}

\begin{figure}
\centering
\includegraphics[width=0.99\textwidth]{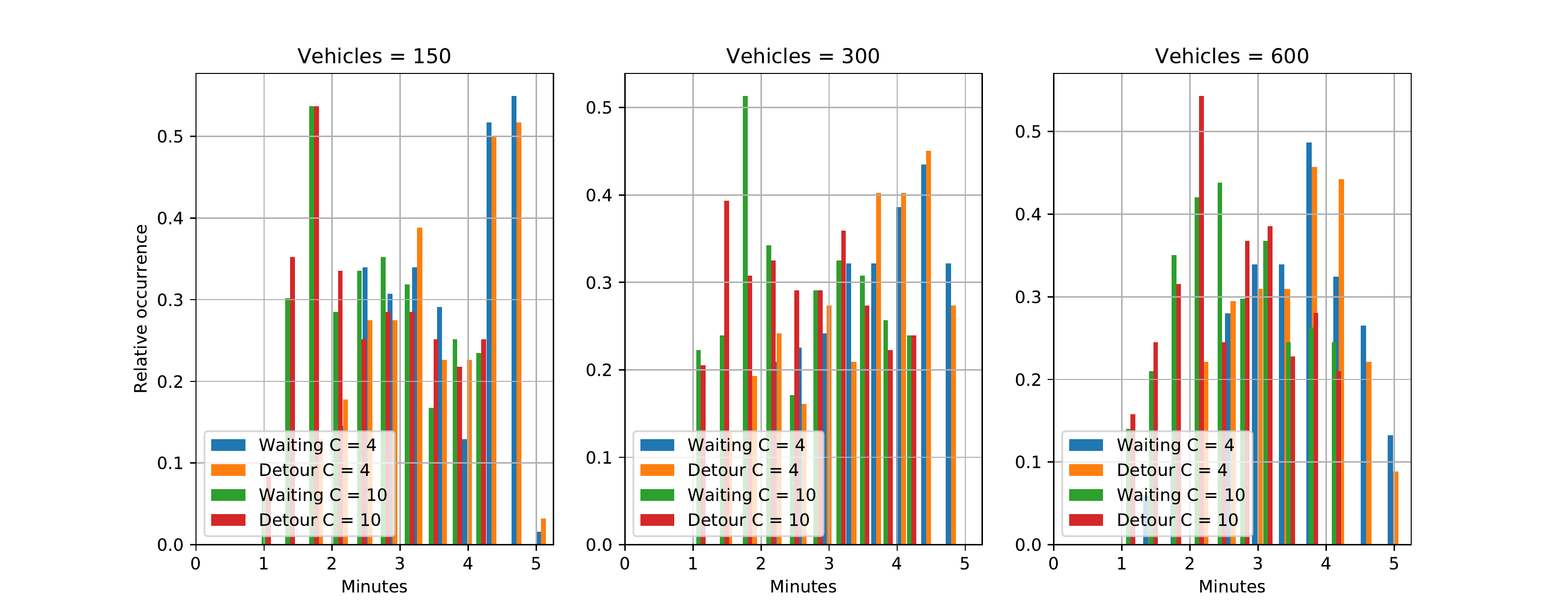}
\caption{Waiting and detour times per customers in the considered time frame (histogram).}
\label{fig.sim-2-bis}
\end{figure}

\begin{figure}
\centering
\includegraphics[width=0.99\textwidth]{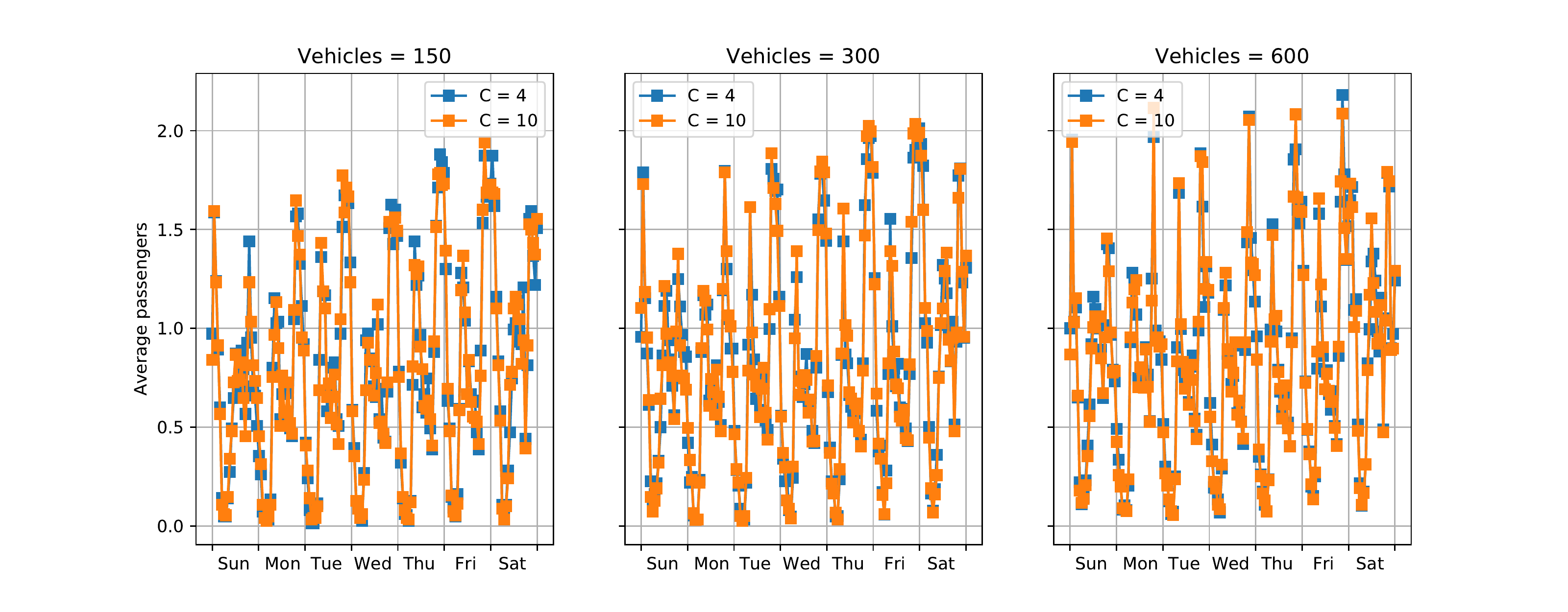}
\caption{Average passengers per vehicle in the considered time frame.}
\label{fig.sim-3}
\end{figure}

\begin{figure}
\centering
\includegraphics[width=0.99\textwidth]{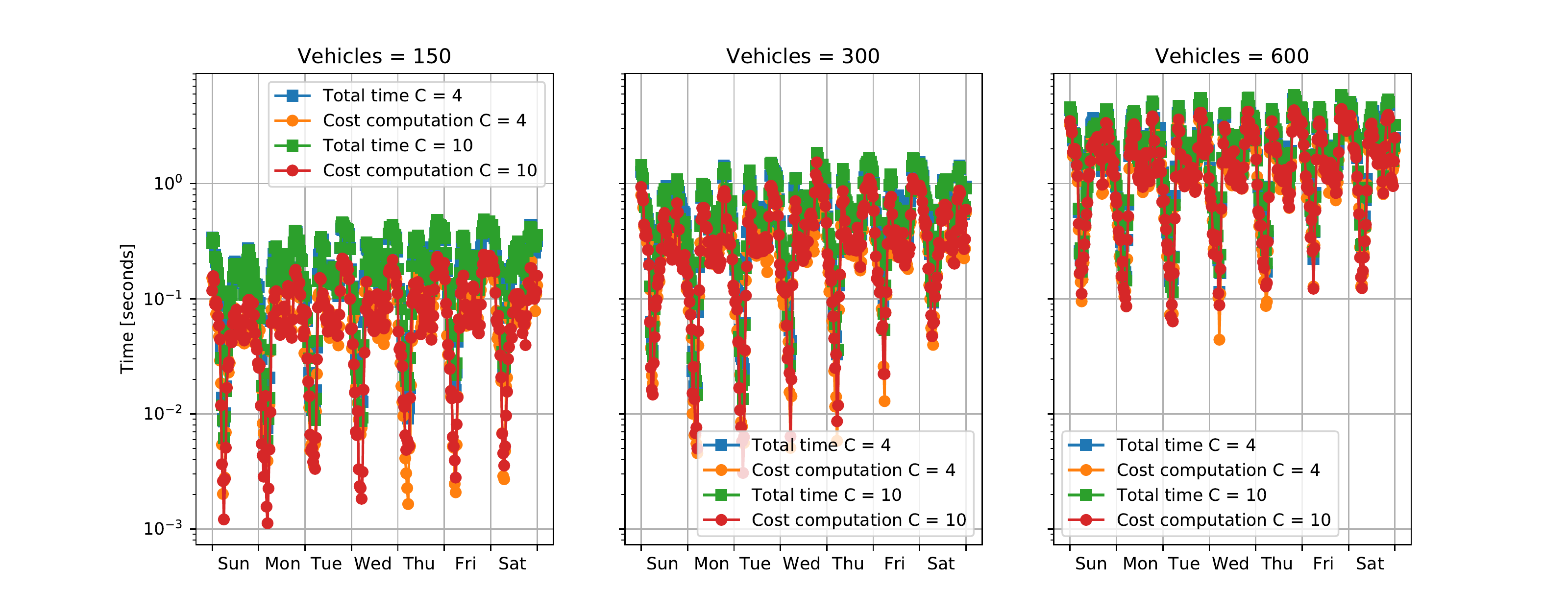}
\caption{Computational time per call of the ridesharing service in the considered time frame. The $y$-axis is in logarithmic scale. }
\label{fig.sim-4}
\end{figure}

\begin{figure}
\centering
\includegraphics[width=0.99\textwidth]{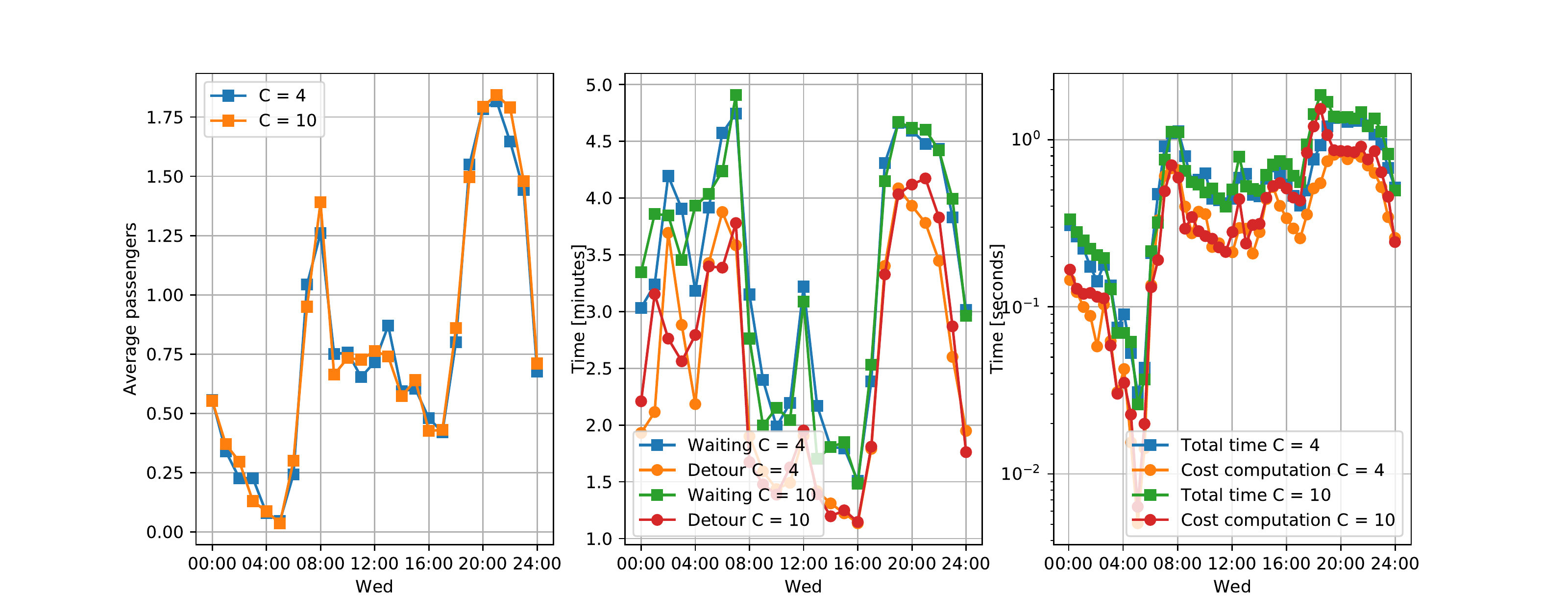}
\caption{Detail for the day of Wednesday for the case of $300$ vehicles.}
\label{fig.sim-100}
\end{figure}

\begin{figure}
\centering
\includegraphics[width=0.99\textwidth]{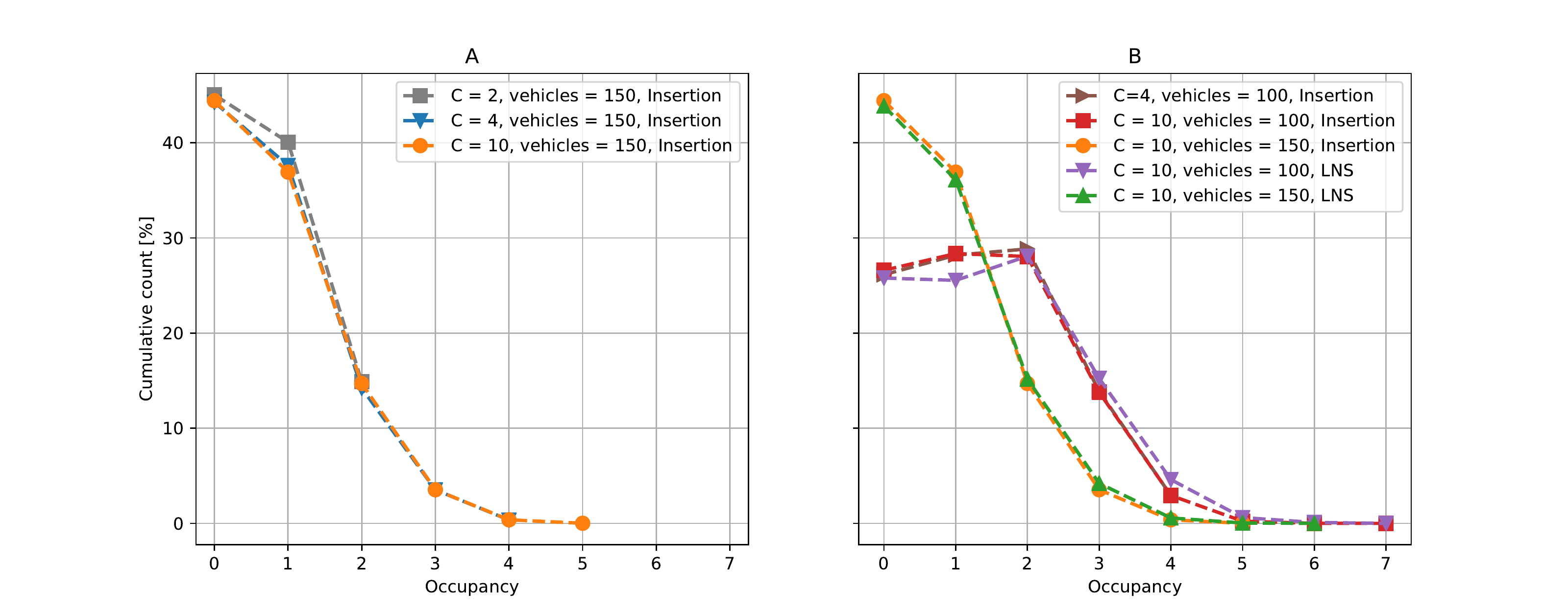}
\caption{Occupancy counts: the $x$-axis represents the number of customers in a vehicle, while the $y$-axis represents the frequency of the occurrence. Figure (A) relates to the case of $150$ vehicles and DARP solved with the insertion heuristic. Figure (B) reports the fleet occupancy for different scenarios with varying capacity, fleet size and DARP algorithm. }
\label{fig.sim-5}
\end{figure}

To assess the algorithm performance, we present:
\begin{enumerate}
\item \revtwo{Numerical simulations with $150$, $300$, and $600$ vehicles, which correspond to $5\%$, $10\%$, and $20\%$ of the maximum fleet size considered in~\cite{Alonso-Mora2017}. For each fleet size, we consider the same percentage of ride requests, see Figures~\ref{fig.sim-1} to \ref{fig.sim-100}, and Tables~\ref{tab:ny-r-1-bis} and \ref{tab:ny-r-1}. These results show that, even for small fleet sizes, if the number of customers is scaled accordingly, then the performance of the service is almost not affected. This advocates the service for small companies and larger ones alike.} 

\item Numerical simulations with variations of the algorithm parameters, namely the maximum number of vehicles considered in the context mapping (max$n$), the cost function, the sampling period, the vehicle capacity, \revtwo{and the distance metric considered in the context mapping module} are reported in Table~\ref{tab:ny-r-2} and  Table~\ref{tab:ny-r-2-bis}. These results show how the different parameters affect the solution quality. 

\item \revtwo{Numerical simulations with time variations of demand, fleet, and congestion are reported in Table~\ref{tab:ny-r-time-v} and Figures~\ref{rain}-\ref{rainstuff}. These results demonstrate the robustness and sensitivity of the algorithm to small and not-so-small time-variations.  
}

\item Numerical simulations comparing the solution yielded by the insertion heuristic and a large neighborhood search algorithm (see Algorithm~\ref{alg:ins-heur} and Algorithm~\ref{alg:lns}), cf. Figure~\ref{fig.sim-5} and Table~\ref{tab:ny-r-3}. This comparison displays how the insertion heuristic is a reasonable and efficient option for our ridesharing solution. 

\item \revtwo{Numerical simulations with the whole ridesharing demand, and a fleet size of $2000$ and $3000$. These results are reported as a more accurate comparison with the ones obtained in~\cite{Alonso-Mora2017}. The findings in Table~\ref{tab:benchmark} show qualitatively similar performance and lower computational time. In particular, we achieve a similar performance as~\cite{Alonso-Mora2017} while our computational time is as much as $4$ times lower than theirs. This is even more impressive, since we implement our algorithm using the Python 2.7 language on a 2.7GHz Intel i5 laptop with 8GB RAM, while in~\cite{Alonso-Mora2017} their algorithm is implemented in C++ on a 24 core machine. 
}  

\end{enumerate}




\paragraph{Basic performance.} 

We present numerical simulations with $150$, $300$, and $600$, which correspond to $5\%$, $10\%$, and $20\%$ of the maximum fleet size considered in~\cite{Alonso-Mora2017}. For each fleet size, we consider the same percentage of ride requests.

In Figure~\ref{fig.sim-1}, we depict the percentage of serviced customers in the case of different vehicle fleet size and vehicle capacity (either $C = 4$, or $C = 10$). The percentage of serviced customers labeled subsequently as service rate (SR) is the percentage of serviced customers with respect to the total number of requests. In this sense, the percentage of refused customers is $100 - \textrm{SR}$. The serviced customers are computed on the time frame from 00:00 hours, Sunday till 23:59 hours, Saturday, to avoid any spurious effect due to the random initialization.  As we can see, in the case of customers accepting the rebalancing vehicles (A), or not (B), the algorithm guarantees a high service level. We notice that one does not need a high service penetration (i.e., market share) to obtain a very high service level: $10\%$ seems more than sufficient. Notice that $300$ vehicles are a little more than the $2\%$ of the current taxi fleet and are enough to handle $10\%$ of the current demand thanks to ridesharing. These results are in line with the ones presented for the whole fleet in~\cite{Alonso-Mora2017} (cf. Table~\ref{tab:ny-r-1}). Since service rates in cases (A) or (B) are similar, the customer preference on rebalancing does not affect considerably the service. This suggests that a combination of demand forecasting techniques and strategies to move idle vehicles towards highly-demanded areas could be beneficial for such ridesharing services. We remark here that the more vehicles we consider (and correspondingly the more customers), the more the service rate increases. This is because requests spatial density increases and vehicles reduce their trip lengths to accommodate the requests. In this sense, the results for $5\%, 10\%, 20\%$ represent lower bounds on the performance of our algorithm in the case of the whole fleet (provided that max$n$ is scaled accordingly).  

Figure~\ref{fig.sim-1}(C) shows the situation in which one does not perform reactive rebalancing of the fleet. In this case, the quality of service is highly deteriorated. In particular, the service rate drops of $25\%$. It is interesting to notice that, for a fleet of $150$ vehicles, the solution for $C=4$ seems better that the one for $C=10$.
\revtwo{This seems counter-intuitive, but may be explained as follows. When $C=4$ and the number of vehicles is relatively limited, it can often happen that the non-idle vehicles are full, and hence the algorithm would ask idle vehicles to reach demanded areas (either in the assignment module or in the reactive rebalancing). On the other hand, when $C=10$ the idle vehicles may not be asked to follow the demand frequently, according to the reactive rebalancing strategy: this behavior might be investigated in a future work with an active rebalancing.}

In Figures~\ref{fig.sim-2} and~\ref{fig.sim-2-bis}, we capture how the constraints on maximum waiting time and detour time are satisfied along the considered time period. In particular, we report the case in which the customers are accepting the rebalancing vehicles (thereby not necessarily satisfying the constraints), and despite this, the average values in hour batches is well below the constraint limits. The total average values are reported in Table~\ref{tab:ny-r-1} and they are similar to the ones reported in~\cite{Alonso-Mora2017}. In Table~\ref{tab:ny-r-1}, ``y'' indicates the strategy of accepting a rebalancing vehicle, while ``n'' the strategy of refusing it. \revtwo{In the histogram plots, we see how the distribution of hourly waiting times and detour times changes for various fleet size. In the case of $C=4$, the distribution peaks shift on the left by increasing the number of vehicles, providing a better service. For $C=10$, the peaks are more on the left than the case $C=4$, so the service is in general better.} \revtwo{This further demonstrates that rebalancing occurs less frequently in the case of limited capacity. It can be observed that as the number of vehicles increases, the distribution gets more spread out.}

Figure~\ref{fig.sim-3} gives an indication of how the fleet is utilized by showing the number of passengers per vehicle in the considered time frame. Given that the average occupancy of the vehicles rarely exceeds $2$, the simple insertion heuristic Algorithm \ref{alg:ins-heur} should give solutions of acceptable quality. This point will be elaborated at the end of the section.

We report in Figure~\ref{fig.sim-4} the computational time spent to run the algorithm. The most computational intensive part is the cost computation, which scales linearly in the number of vehicles and in ${\max}n$. All the other parts of the algorithm are significantly faster, even the linear assignment problem, which takes less than $5$ seconds to solve problems of size $3000$ in our implementation. As we see, vehicle fleets up to size $600$ can be implemented in real-time on a simple laptop (rem: the batch period is $10$ seconds). 


\rev{
Figure~\ref{fig.sim-100} represents a zoom in for the day of Wednesday, so to better appreciate hourly variations. Here we see that the results are qualitatively similar between a capacity of $C=4$ and $C=10$. We further appreciate the hourly variations, in particular, waiting and detour time are larger in the morning commute and in the evening when people come home from work or go out. Occupancy is related to the number of requests, which are at their lowest value around 04:00 hours and at their peak around 20:00 hours. 

In Table~\ref{tab:ny-r-1-bis}, we report heat maps for the occupancy of vehicles in a $300$, $C=4$ vehicle scenario for the considered week. We can appreciate daily and hourly variations. In particular, in the morning commute the vehicles are picking up customers in upper Manhattan, while delivering them in later hours in mid and lower Manhattan. During the day the busiest areas are still in the mid and lower parts. 

}

\begin{table}
\caption{Heat maps for the occupancy of vehicles in a $300$, $c=4$ vehicle scenario for the considered week.}
\label{tab:ny-r-1-bis}
 \sffamily
 {\renewcommand{\arraystretch}{1.05}
\begin{tabular}{p{.01\textwidth}p{.0975\textwidth}p{.0975\textwidth}p{.0975\textwidth}p{.0975\textwidth}p{.0975\textwidth}p{.0975\textwidth}p{.0975\textwidth}p{.0975\textwidth}}
 &\centering 01:00 & \centering 04:00 & \centering 07:00 & \centering 10:00 & \centering 13:00 & \centering 16:00 & \centering 19:00 &  \hskip0.5cm 22:00  \\ \toprule
{\rotatebox[origin=c]{90}{\hskip.75cm Monday \hskip1.5cm Sunday\hskip.75cm}} & \multirow{8}{5cm}{\vskip-3cm\includegraphics[width = 0.96\textwidth]{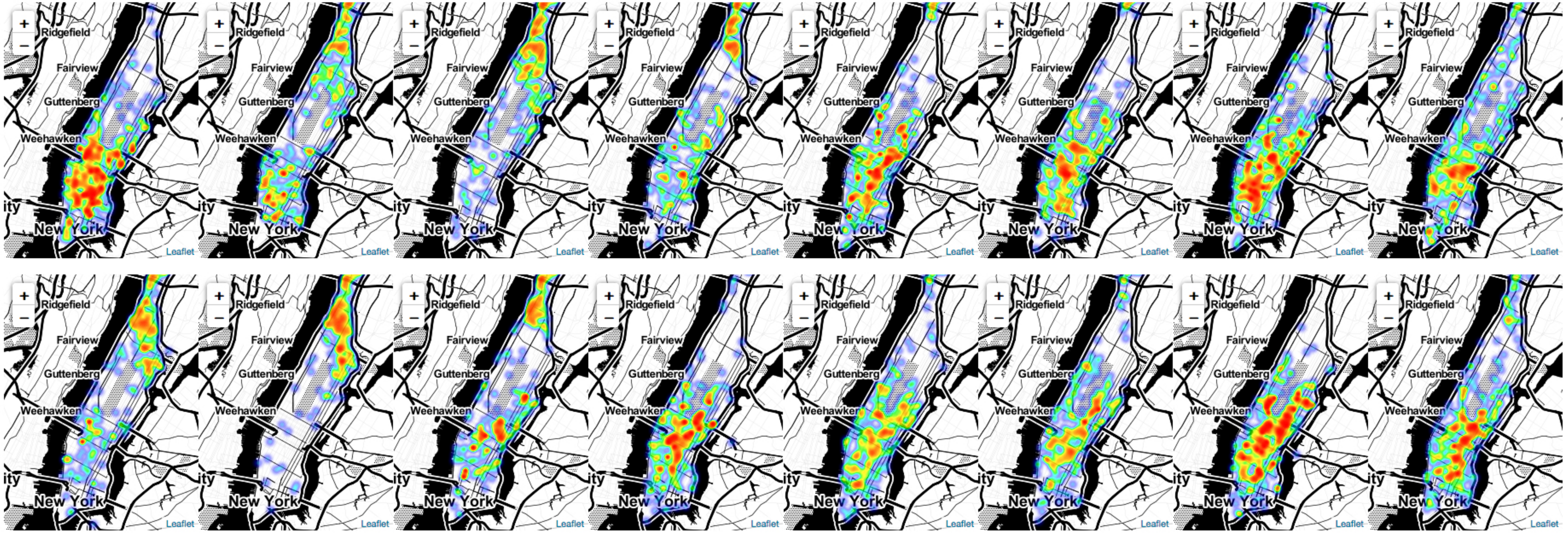}}  \\
{\rotatebox[origin=c]{90}{\hskip.75cm Wednesday \hskip1.3cm Tuesday\hskip.95cm}} & \multirow{8}{5cm}{\vskip-3cm\includegraphics[width = 0.96\textwidth]{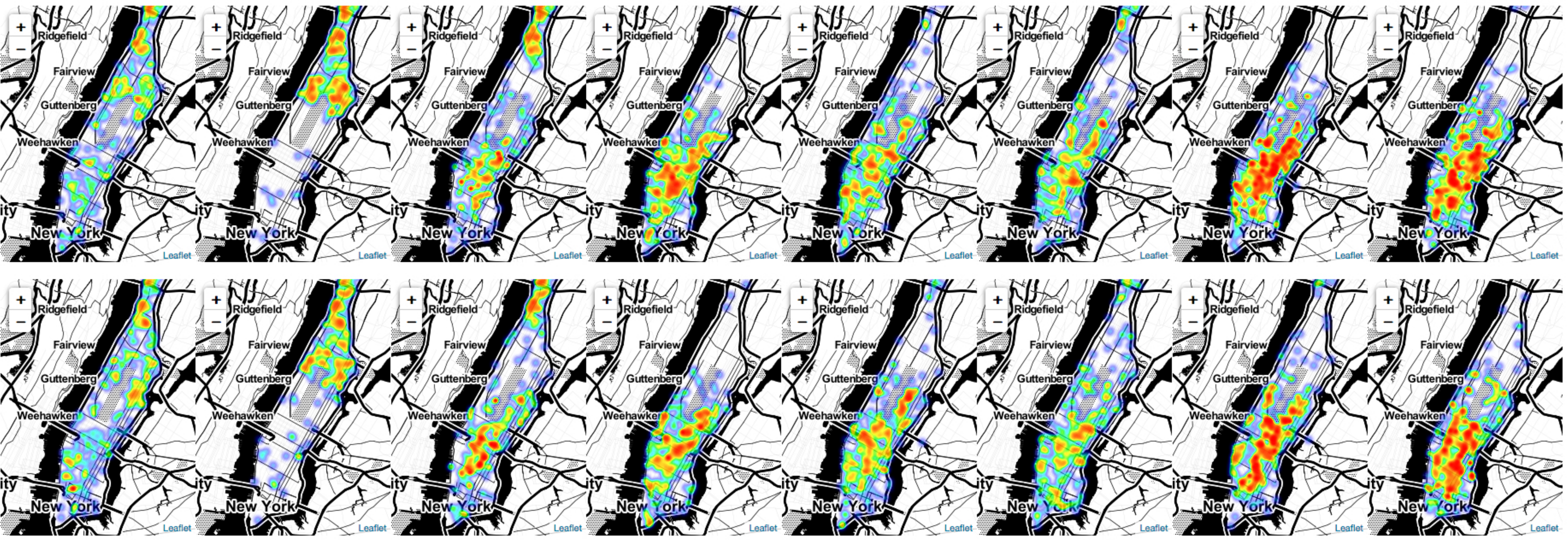}}  \\
{\rotatebox[origin=c]{90}{\hskip.75cm Friday \hskip1.5cm Thursday\hskip.75cm}} & \multirow{8}{5cm}{\vskip-3cm\includegraphics[width = 0.96\textwidth]{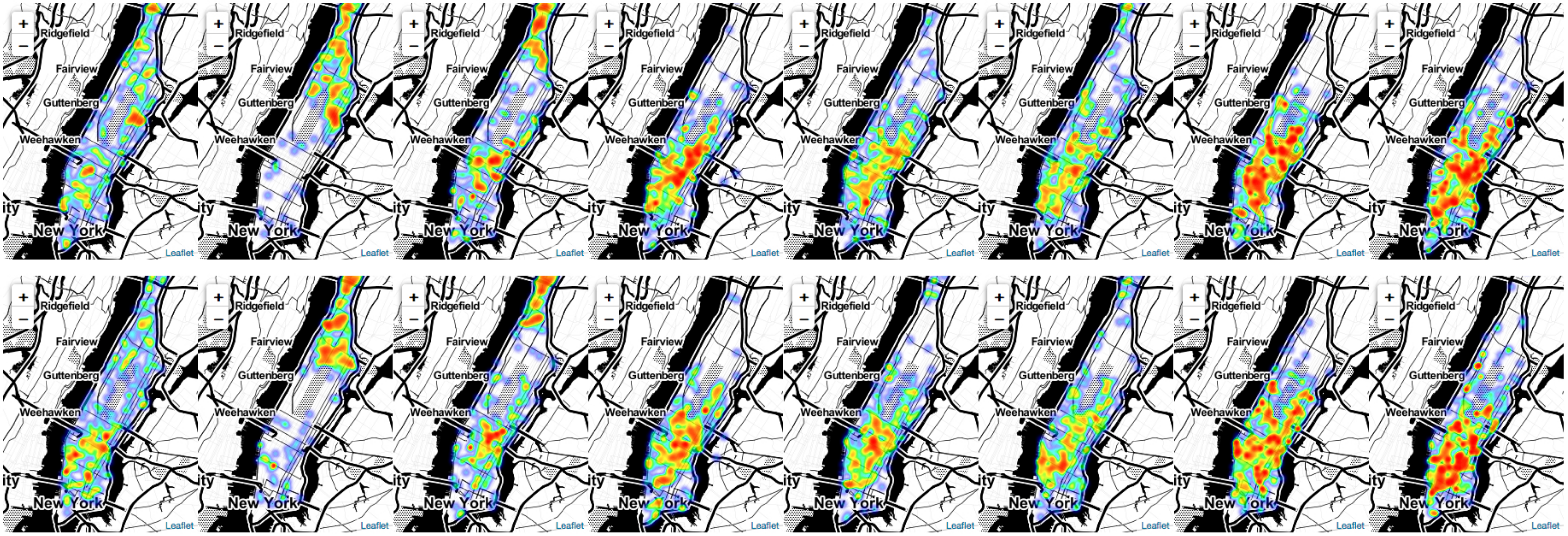}}  \\
{\rotatebox[origin=c]{90}{\hskip.75cm Saturday \hskip.75cm}} & \multirow{8}{5cm}{\vskip-3cm\includegraphics[width = 0.96\textwidth]{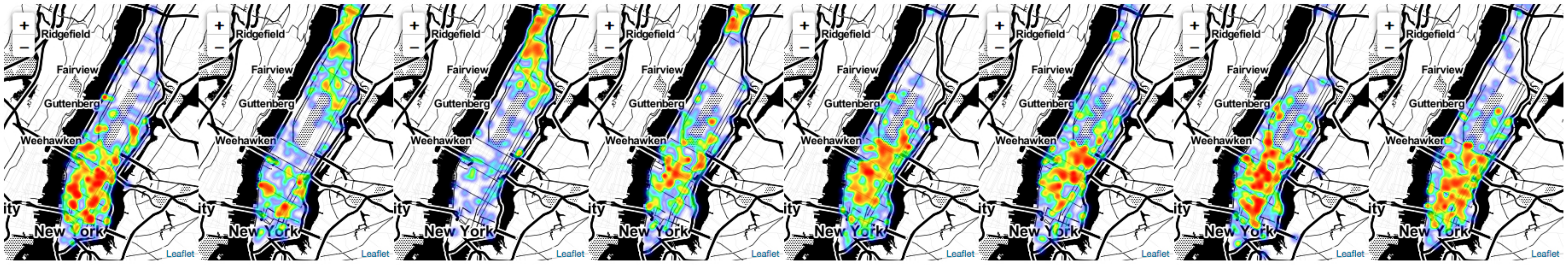}}  
\end{tabular} } 
\end{table}

\rev{

Table~\ref{tab:ny-r-1} summarizes the results and the comparison with the results in~\cite{Alonso-Mora2017}. The columns of the table are the number of vehicles, the percentage of the vehicles w.r.t. the $3000$ vehicles used in~\cite{Alonso-Mora2017}, the percentage of requests considered w.r.t. the total amount of requests, the capacity $C$, the maximum number of considered vehicles per requests (max$n$), the cost function (rem: TD stands for total trip duration), the sampling period, the service rate (SR), the average waiting time with (y) and without (n) accepting rebalanced vehicles, the average detour time with (y) and without (n) accepting rebalanced vehicles, and the computational time. 

Note that the cost function of~\cite{Alonso-Mora2017} is indicated as $C(\Sigma)$ as in the cited paper and it is a combination of detour and waiting time and it adds a penalty for refusing customers. Note also that in~\cite{Alonso-Mora2017} customers never accept rebalanced vehicles. 

As we see, the performance in all the cases is very similar and we also compare similarly to~\cite{Alonso-Mora2017}. We note here that, service rate (SR) from $C=4$ and $C=10$ is similar to~\cite{Alonso-Mora2017}, and increasing, as expected, with the number of vehicles. SR is increasing since the requests are spatially denser and the trip length is shorter, therefore the vehicles are more able to handle new requests. On the other hand, since we optimize for trip duration (TD), we naturally push for low occupancy and therefore $C=4$ has similar performance than $C=10$. Actually, in some cases $C=4$ performs just slightly better than $C=10$, which is somehow counter-intuitive, but can happen (especially with ``few'' vehicles and when $C=4$ is ``sufficient''), since in the low capacity case more vehicles are being moved and rebalancing happens less often. This is not the case for a number of vehicle of $100$ with the same number of requests (see Table~\ref{tab:ny-r-3}). 


\begin{table}
\footnotesize
\centering
\caption{\revtwo{Results for different values of vehicles, customers, and capacity and comparison with~\cite{Alonso-Mora2017} (indicated by $^*$).}}
\label{tab:ny-r-1}
\begin{tabular}{ccccccccccccc}
\toprule
vehicles & vehicles & customers & $c$ & max$n$ & cost & $h$ & SR & waiting & waiting  &  detour  & detour & comp. time  \\ 
 & [\%] & [\%] & &  & &[s] & [\%] & y [min] & n [min]  &  y [min] & n [min] & [s]  \\ 
\toprule
\rowcolor{Gray}
150 & 5.0 & 5.0 & 4 & 8 & TD & 10 & 95.75 & 3.54 & 3.43 & 2.53 & 2.53 & 0.17 \\  
150 & 5.0 & 5.0 & 10 & 8 & TD & 10 & 95.68 & 3.58 & 3.46 & 2.50 & 2.50 & 0.18 \\ 
\rowcolor{Gray}
300 & 10.0 & 10.0 & 4 & 16 & TD & 10 & 98.88 & 3.46 & 3.38 & 2.60 & 2.60 & 0.56 \\  
300 & 10.0 & 10.0 & 10 & 16 & TD & 10 & 98.87 & 3.45 & 3.37 & 2.60 & 2.59 & 0.58 \\ 
\rowcolor{Gray}
600 & 20.0 & 20.0 & 4 & 32 & TD & 10 & 99.86 & 3.43 & 3.37 & 2.65 & 2.65 & 2.05 \\  
600 & 20.0 & 20.0 & 10 & 32 & TD & 10 & 99.88 & 3.39 & 3.33 & 2.63 & 2.63 & 2.11 \\ 
\midrule
$^*$, 3000 & 100.0 & 100.0 & 2 & - & $C(\Sigma)$ & 30 & 94.21 & - & 3.19 & - & 1.46 & 31.38 \\ 
\rowcolor{Gray} 
$^*$, 3000 & 100.0 & 100.0 & 4 & - & $C(\Sigma)$ & 30 & 97.91 & - & 2.70 & - & 2.28 & 51.55 \\  
$^*$, 3000 & 100.0 & 100.0 & 10 & - & $C(\Sigma)$ & 30 & 98.58 & - & 2.56 & - & 2.74 & 60.39 \\  
\toprule
\end{tabular}
\end{table}

}

\rev{
\paragraph{Performance with different choice of parameters.} We run simulations varying a number of significant algorithm parameters and we obtain the results reported in Tables \ref{tab:ny-r-2} and \revtwo{\ref{tab:ny-r-2-bis}}. In particular, keeping the number of vehicles and requests fixed, we make the following experiments.

\begin{enumerate}
\item The maximum number of vehicles considered per request (max$n$) is varied in the set $6,8,10,12$. In this case, the service rate (SR) improves, while computational time increases, as expected. The choice of max$n$ of $8$ is motivated by trading-off performance with computational time. 

\item We tested the cost functions TD (total trip duration), WT (total waiting time for customers), and DT (total detour time for vehicles). Formally, WT is the total waiting time for customers of a given solution, while DT is the difference between the total trip time before the assignment and after. As we see, TD is the most performing cost function, keeping trip lengths small and making vehicles more available. WT trades-off short waiting times for long detour times, while DT does the opposite. Note that, even though WT minimizes waiting time at every time, given the myopic strategy we are using, there is no guarantee that the waiting times are minimized on average over a week time; actually in this case, TD does a better job. The fact that the computational times are larger for WT and DT implies that in these cases, the vehicles have longer trip times, so they need to keep track of more customers; while naturally TD does not. 

\item The sampling period is varied in the set $5,10,20,30$ seconds. We see that a sampling period of $5$~s offers the best service rate (SR), because it is able to assign more customers to the same vehicle over multiple time steps (indicated by larger detour times). However, larger sampling periods, e.g., $30$~s offers lower detour times with comparable SR. The $30$~s one has an increased SR, since within a larger time frame more requests are available and one can find better request-to-vehicle assignments. The fact that the SR is almost constant while the other indicators change confirms the similar results obtained in~\cite{Alonso-Mora2017}. 

\item The capacity is varied in the set $2,4,10$. As we see, performance increase with the capacity to peak at $C=4$. As discussed in the previous subsection, the case $C=10$ is special for this scenario with a low number of vehicles. In particular, as said, service rate is very similar, since our cost function TD favors small trip durations. This can be even better seen in Figure~\ref{fig.sim-5}(A), where we plot the frequencies of average occupancies in the instance with $150$ vehicles. Since there are no more than $5$ customers on a vehicle route, the solutions obtained for $C=4$ and $C=10$ turn out to be very similar. When the number of vehicles is $100$, Figure~\ref{fig.sim-5}(B) show that the vehicle occupancy can be of $7$ customers. This motivated us to test a more advanced DARP heuristic, i.e., the LNS, and we discuss the performance of DARP algorithms in the next subsection.

\revtwo{
\item The distance metric in the context mapping module is varied to Manhattan distance and to shortest-path distance (the baseline is Euclidean distance). As we see, the difference in service rate and timing is very small. It appears that the Manhattan distance obtains the best service rate, although it degrades the waiting and detour times. Shortest-path seems to perform the worst, once again advocating that the best thing that one can do at a given time, may be not the best in the long run. }
\end{enumerate}
}

\begin{table}
\footnotesize
\centering
\caption{Results for different values of max$n$, cost function, sampling period, capacity.}
\label{tab:ny-r-2}
\begin{tabular}{ccccccccccccc}
\toprule
vehicles & vehicles & customers & $c$ & max$n$ & cost & $h$ & SR & waiting & waiting  &  detour  & detour & comp. time  \\ 
 & [\%] & [\%] & &  & &[s] & [\%] & y [min] & n [min]  &  y [min] & n [min] & [s]  \\ 
\toprule
150 & 5.0 & 5.0 & 4 & 6 & TD & 10 & 94.29 & 3.63 & 3.48 & 2.59 & 2.57 & 0.14 \\ 
\rowcolor{Gray}
150 & 5.0 & 5.0 & 4 & 8 & TD & 10 & 95.75 & 3.54 & 3.43 & 2.53 & 2.53 & 0.17 \\ 
150 & 5.0 & 5.0 & 4 & 10 & TD & 10 & 96.82 & 3.52 & 3.43 & 2.50 & 2.50 & 0.20 \\ 
150 & 5.0 & 5.0 & 4 & 12 & TD & 10 & 97.53 & 3.52 & 3.44 & 2.54 & 2.54 & 0.23 \\
\midrule
\rowcolor{Gray}
150 & 5.0 & 5.0 & 4 & 8 & TD & 10 & 95.75 & 3.54 & 3.43 & 2.53 & 2.53 & 0.17 \\
150 & 5.0 & 5.0 & 4 & 8 & WT & 10 & 92.67 & 3.89 & 3.82 & 4.06 & 4.08 & 0.21 \\ 
150 & 5.0 & 5.0 & 4 & 8 & DT & 10 & 95.57 & 4.45 & 4.32 & 3.85 & 3.85 & 0.26 \\  \midrule
150 & 5.0 & 5.0 & 4 & 8 & TD & 5 & 96.12 & 3.51 & 3.40 & 3.85 & 3.87 & 0.10 \\ 
\rowcolor{Gray}
150 & 5.0 & 5.0 & 4 & 8 & TD & 10 & 95.75 & 3.54 & 3.43 & 2.53 & 2.53 & 0.17 \\
150 & 5.0 & 5.0 & 4 & 8 & TD & 20 & 95.71 & 3.62 & 3.52 & 2.63 & 2.62 & 0.32 \\ 
150 & 5.0 & 5.0 & 4 & 8 & TD & 30 & 95.96 & 3.63 & 3.52 & 2.58 & 2.59 & 0.39 \\ 
\midrule
150 & 5.0 & 5.0 & 2 & 8 & TD & 10 & 93.82 & 3.65 & 3.54 & 3.31 & 3.32 & 0.12 \\ 
\rowcolor{Gray}
150 & 5.0 & 5.0 & 4 & 8 & TD & 10 & 95.75 & 3.54 & 3.43 & 2.53 & 2.53 & 0.17 \\
150 & 5.0 & 5.0 & 10 & 8 & TD & 10 & 95.62 & 3.62 & 3.51 & 2.62 & 2.62 & 0.18 \\ 
\toprule
\end{tabular}
\end{table}

\begin{table}
\footnotesize
\centering
\caption{\revtwo{Results for different distance metrics in the context mapping module, max$n$ = 8, cost = TD. }}
\label{tab:ny-r-2-bis}
\begin{tabular}{cccccccccccc}
\toprule
vehicles & vehicles & customers & $c$  & distance & $h$ & SR & waiting & waiting  &  detour  & detour & comp. time  \\ 
 & [\%] & [\%] & &   &[s] & [\%] & y [min] & n [min]  &  y [min] & n [min] & [s]  \\ 
\toprule
\rowcolor{Gray}
150 & 5.0 & 5.0 & 4  & Euclidean & 10 & 95.75 & 3.54 & 3.43 & 2.53 & 2.53 & 0.17 \\ 
150 & 5.0 & 5.0 & 4  & Manhattan & 10 & 95.77 & 3.61 & 3.49 & 2.55 & 2.55 & 0.17 \\ 
150 & 5.0 & 5.0 & 4  & Shortest-path & 10 & 95.70 & 3.59 & 3.48 & 2.59 & 2.59 & 0.19 \\
\toprule
\end{tabular}
\end{table}

\begin{table}
\footnotesize
\centering
\caption{Difference between insertion heuristic and large neighborhood search, cost = TD.}
\label{tab:ny-r-3}
\begin{tabular}{ccccccccccccc}
\toprule
vehicles & vehicles & customers & $c$ & max$n$ & method & $h$ & SR & waiting & waiting  &  detour  & detour & comp. time  \\ 
 & [\%] & [\%] & &  & &[s] & [\%] & y [min] & n [min]  &  y [min] & n [min] & [s]  \\ 
\toprule
\rowcolor{Gray}
150 & 5.0 & 5.0 & 10 & 8 & Ins & 10 & 95.75 & 3.54 & 3.43 & 2.53 & 2.53 & 0.17 \\ 
150 & 5.0 & 5.0 & 10 & 8 & LNS & 10 & 95.70 & 3.57 & 3.47 & 2.77 & 2.77 & 0.49 \\ 
\midrule
100 & 3.3 & 5.0 & 4 & 8 & Ins & 10 & 75.16 & 4.40 & 4.29 & 3.80 & 3.80 & 0.27 \\
\rowcolor{Gray}
100 & 3.3 & 5.0 & 10 & 8 & Ins & 10 & 75.27 & 4.37 & 4.27 & 3.80 & 3.80 & 0.28\\ 
100 & 3.3 & 5.0 & 10 & 8 & LNS & 10 & 75.37 & 4.29 & 4.20 & 4.12 & 4.12 & 0.57 \\  
\toprule
\end{tabular}
\end{table}

\rev{
\paragraph{Performance with different DARP algorithms.} 
We further test our ridesharing algorithm varying the DARP solution approach, and report the results in Table~\ref{tab:ny-r-3}. For instances with $150$ vehicles and $C=10$, the LNS deteriorates the SR of a percentage of 0.05\%. As shown in Figure~\ref{fig.sim-5}(A), the occupancy of vehicle seldom exceeds $4$, therefore LNS and insertion heuristic solutions are often the same. However, it can happen that the myopic solution of LNS yields a slightly worse result that those of the insertion heuristic in the overall one-week service.

For both cases $n=100$ and $n=150$, the LNS requires larger detour times for vehicles, in order to select routes with a shorter duration. The computational time of the LNS is around the double of the insertion heuristic. On the instance with $100$ vehicles, the LNS solution provides a slightly better service rate and smaller waiting times for customers: as shown in Figure~\ref{fig.sim-5}(B), in instances with a smaller fleet the average occupancy increases, therefore the quality of LNS becomes more evident. It can also be observed that, with respect to the insertion heuristic, the LNS tends to favor more routes with an occupancy larger than $2$. To sum up, the insertion heuristic proves to be the best DARP algorithm, as it finds an acceptable trade-off between level of service and computational time. This is expected in our instances, because of the low number of scheduled requests on the vehicle routes.
}

%
%

\revtwo{
\paragraph{Sensitivity to time variations.}

We move on to test the sensitivity of the algorithm to time variations. In particular, we present tests for:

\begin{enumerate}
\item Time variations on the demand. We obtained the precipitation profile for the tested week in Mahnattan from Wolfram Alpha, which is also reported in Figure~\ref{rain}. We increased the amount of requests proportionally to the precipitation amount minute per minute, so that the maximum precipitation amount on Wednesday causes an increase in demand of either $20\%$ (Simulation 1) or $40\%$ (Simulation 2). Table~\ref{tab:ny-r-time-v} and Figure~\ref{rainstuff} report our findings.  As expected, the performance of the ridesharing algorithm is negatively affected by the increased demand. However, if the increase of demand occurs only during a limited time, the decrease in performance w.r.t. the baseline of no rain is very small, on average. A zoom-in analysis in Figure~\ref{rainstuff} reveals that the most affected performance metric is the number of refused customers, while waiting times are only partially affected. In particular, around 12:00 the total percentage of refused customers goes from $0.006\%$ (baseline) to $0.065\%$ in Simulation 1 (which corresponds to refusing around $8\%$ of incoming requests, instead than less than $1\%$ in the baseline) and to $0.185\%$ in Simulation 2 (which means refusing $19\%$ of incoming requests). These statistics quantify the flexibility of the algorithm to time-varying demand, and it is interesting to note that the service level remains above $95\%$.

\item Time variations on the taxi fleet size. We varied the taxi fleet online. In particular, every 30 minutes we randomly set $20\%$ (Simulation 1) or $40\%$ (Simulation 2) of the idle vehicles as taking a 30 minutes break. The simulations accounts for fluctuations in fleet size due to driver shifts, which are realistic from a business perspective. More vehicles are taking a break when the number of requests is low (because more vehicles are idle), than when the demand is high. Table~\ref{tab:ny-r-time-v} reports our findings. The performance of the ridesharing algorithm degrades, but not considerably, which is encouraging for real systems where drivers might take breaks. As one can also observe, waiting and detour times increase (as well as the computational time), which indicates that vehicles try to fit more passengers per ride.  

\item Time variations for congestion. We obtained the average velocity in Manhattan hour-by-hour by processing the Taxi data. In particular, for each trip, we divided the shortest path distance with the real duration of the trip. Averaging across all the trips in a given hour, we extract the average hourly velocity. We computed the average over the week and looked at the ratio between hour velocity and average one: during the night, one can obtain ratio as high as $1.6$ (hour/average), while in busy time one can go as low as $0.6$. 

We multiplied the travel times that OSRM computes by the hourly velocity ratio, so as to simulate time-varying traffic conditions. We then consider two sub-cases: \emph{(i)} the ridesharing system knows/learns the hourly ratios and act accordingly, \emph{(ii)} the ridesharing system does not know the hourly ratios and schedules based on the OSRM travel times. The results can be found in Table~\ref{tab:ny-r-time-v}. Comparing \emph{(i)} to \emph{(ii)}, the performance decreases, since high congestion goes hand in hand with high demand; performance is better for larger fleet sizes (which are then more robust). Case \emph{(ii)} is more taxing than case \emph{(i)} in terms of waiting and detour times, since in the latter one knows what to expect and plans accordingly.

\end{enumerate}

\begin{figure}
\centering
\includegraphics[width = 0.99\textwidth]{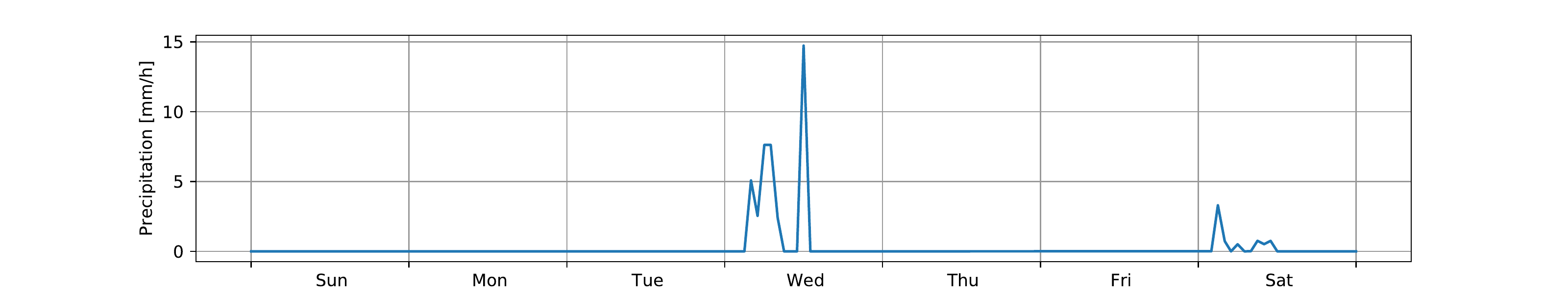}
\caption{\revtwo{Precipitation in New York City during the studied week.}}
\label{rain}
\end{figure}

\begin{figure}
\centering
\includegraphics[width = 0.99\textwidth]{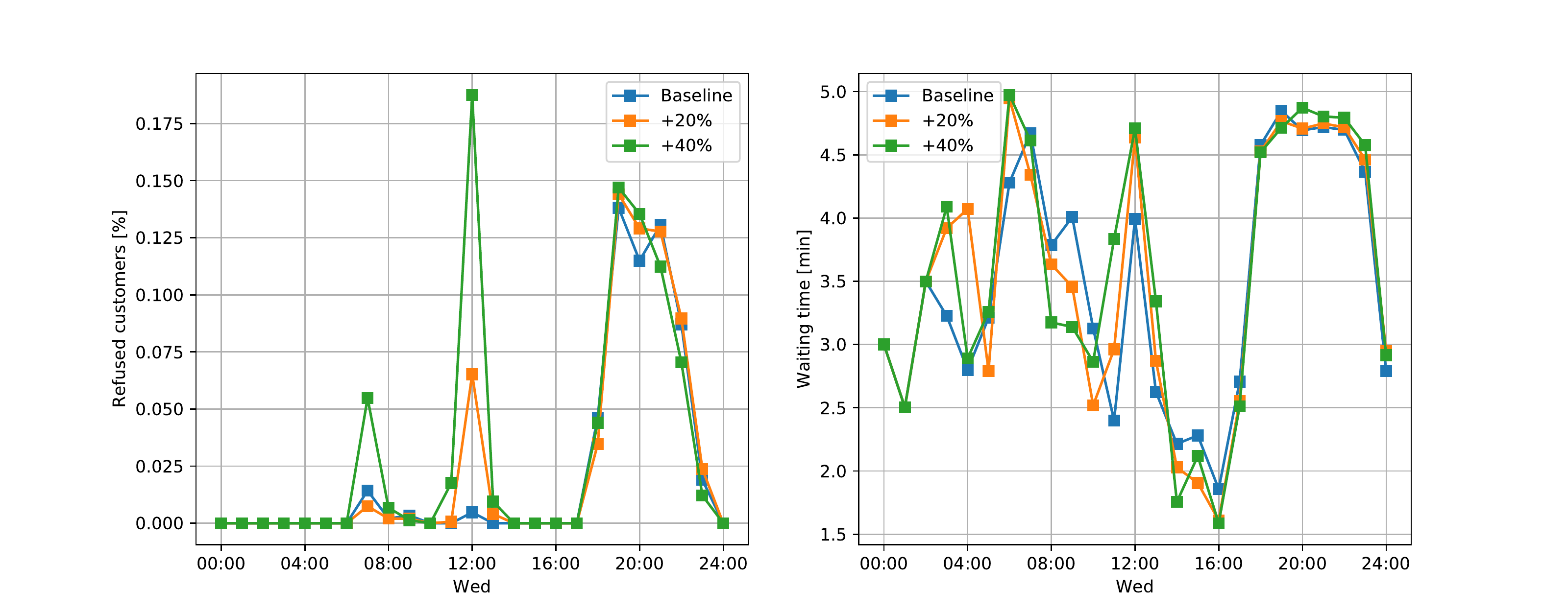}
\caption{\revtwo{Refused customers and waiting time for the considered time-varying demand scenarios. }}
\label{rainstuff}
\end{figure}

\begin{table}
\footnotesize
\centering
\caption{\revtwo{Time variations effects. $c = 4$, $h = 10$~s, cost = TD, max$n$ is 8 for 150 vehicles, 16 for 300 vehicles. }}
\label{tab:ny-r-time-v}
\begin{tabular}{cccccccccc}
\toprule
vehicles & vehicles & customers  & time variations  & SR & waiting & waiting  &  detour  & detour & comp. time  \\ 
 & [\%] & [\%] & & [\%] & y [min] & n [min]  &  y [min] & n [min] & [s]  \\ 
\toprule
\rowcolor{Gray}
150 & 5.0 & 5.0 & None & 95.75 & 3.54 & 3.43 & 2.53 & 2.53 & 0.17 \\ 
\rowcolor{Gray}
300 & 10.0 & 10.0 & None &  98.88 & 3.46 & 3.38 & 2.60 & 2.60 & 0.56 \\ 
\midrule
150 & 5.0 & 5.0 & Demand (+20\%)  & 95.72 & 3.61 & 3.49 & 2.56 & 2.56& 0.17 \\ 
150 & 5.0 & 5.0 & Demand (+40\%)  & 95.57 & 3.60 & 3.49 & 2.58 & 2.58& 0.18 \\ 
\midrule
150 & 5.0 & 5.0 & Fleet-size (-20\%) & 95.66 & 3.61 & 3.51 & 2.60 & 2.60 & 0.18 \\
150 & 5.0 & 5.0 & Fleet-size (-40\%) & 95.33 & 3.70 & 3.59 & 2.68 & 2.68 & 0.25 \\ 
\midrule
150 & 5.0 & 5.0 & Congestion \emph{(i)} & 89.62 & 4.10 & 3.86 & 3.01 & 3.01 & 0.19 \\ 
300 & 10.0 & 10.0 & Congestion \emph{(i)} & 94.24 & 3.98 & 3.84 & 3.09 & 3.09 & 0.63 \\ 
150 & 5.0 & 5.0 & Congestion \emph{(ii)} & 87.91 & 4.50 & 4.28 & 3.76 & 3.76 & 0.20 \\ 
300 & 10.0 & 10.0 & Congestion \emph{(ii)} & 91.91 & 4.41 & 4.27 & 3.87& 3.87 & 0.67 \\ 
\toprule
\end{tabular}
\end{table}

\paragraph{Results with the entire fleet and demand.} We present here the results for our algorithm with the entire ridesharing fleet and demand. Results are collected in Table~\ref{tab:benchmark}. As we can see, the results obtained with our algorithm are in line with the ones in~\cite{Alonso-Mora2017} in terms of service rate, but with a significant reduction of computational times. In particular, in the case of 3000 vehicles and using the same sampling period of $30$~s, we obtain a better service rate and a computational time reduction of $4$ times, which is compatible with the real-time sampling frequency (i.e., average computational time lower than the sampling period). 
 The result is remarkable, since our algorithm is implemented in Python 2.7 language on a 2.7GHz Intel i5 laptop with 8GB RAM, while the algorithm of~\cite{Alonso-Mora2017} is implemented in C++ and runs on a 24 core machine.

As we can see from the results, varying $\max$n, we can achieve near real-time servicing with our algorithms in all cases considered, with more than reasonable service rate and waiting times. Increasing $\max$n, the service rate improves, and naturally the computational times increases.

Finally, we remark that the service rate for our approach, for $300$ vehicles, is consistently better than that of~\cite{Alonso-Mora2017}. We believe this might be caused by two effects: (a) the different routing engine used (b) the choice of the cost function. While it is hard to have visibility on the first point, we may argue on the second. The cost $C(\Sigma)$ used in~\cite{Alonso-Mora2017} penalizes refusing customers in a myopic way, while we do not. This gives $C(\Sigma)$ a better service rate at a given time, but saying nothing about the full horizon. That is: not refusing any customers now while further jamming already packed vehicles, may be detrimental in the long term. On the other hand, it appears that a simpler cost, like TD, knows better when to refuse customers, so that on average, the service rate is higher. We think, for dynamic strategies, \emph{knowing when to refuse} is key in obtaining good solutions.

\paragraph{Remark.} We note here that the method in~\cite{Alonso-Mora2017} speed depends upon the different thresholds and time-cuts that they use. It could be made faster with very tight constraints, or run at optimality if enough computational resources and time were available. In the above comparison, we have taken the results that the authors provide in their paper.  Nonetheless, what it is important to note here is that it seems that a more elaborated assignment with time cut-offs (their method) does not necessarily perform better than a (basically) greedy assignment (our method) when recomputed at high rate. 
}

\begin{table}
\footnotesize
\centering
\caption{\revtwo{Results for the entire demand and comparison with~\cite{Alonso-Mora2017} (indicated by $^*$).}}
\label{tab:benchmark}
\begin{tabular}{ccccccccccccc}
\toprule
vehicles & vehicles & customers & $c$ & max$n$ & cost & $h$ & SR & waiting & waiting  &  detour  & detour & comp. time  \\ 
 & [\%] & [\%] & &  & &[s] & [\%] & y [min] & n [min]  &  y [min] & n [min] & [s]  \\ 
\toprule
\rowcolor{Gray}  
\revtwo{2000} & 66.7 & 100.0 & 4 & 25 & TD & 10 & 92.10 &  3.95 & 3.88 & 3.41 & 3.40 & 10.10 \\
\revtwo{2700} & 90.0 & 100.0 & 4 & 25 & TD & 10 & 99.51 & 3.35 & 3.27 & 2.67 & 2.67 & 7.96 \\
\rowcolor{Gray} 
\revtwo{3000} & 100.0 & 100.0 & 4 & 8 & TD & 10 & 98.65 & 3.11 & 3.10 & 2.39 & 2.39 & 4.25 \\  
\revtwo{3000} & 100.0 & 100.0 & 4 & 25 & TD & 10 & 99.87 & 3.31 & 3.23 & 2.60 & 2.59 & 7.87 \\ 
\rowcolor{Gray}
\revtwo{3000} & 100.0 & 100.0 & 4 & 75 & TD & 10 & 99.99 & 3.23 & 3.16 & 2.75 & 2.74 & 19.80 \\
\revtwo{3000} & 100.0 & 100.0 & 4 & 8 & TD & 30 & 99.09 & 3.14 & 3.04 & 2.48 & 2.47 & 12.77 \\  
\toprule
\rowcolor{Gray} 
$^*$, 2000 & 66.7 & 100.0 & 4 & - & $C(\Sigma)$ & 30 & 93.70 & - & 3.28 & - & 3.29 & 57.55 \\ 
$^*$, 3000 & 100.0 & 100.0 & 2 & - & $C(\Sigma)$ & 30 & 94.21 & - & 3.19 & - & 1.46 & 31.38 \\ 
\rowcolor{Gray} 
$^*$, 3000 & 100.0 & 100.0 & 4 & - & $C(\Sigma)$ & 30 & 97.91 & - & 2.70 & - & 2.28 & 51.55 \\  
$^*$, 3000 & 100.0 & 100.0 & 10 & - & $C(\Sigma)$ & 30 & 98.58 & - & 2.56 & - & 2.74 & 60.39 \\  
\toprule
\end{tabular}
\end{table}

\subsection{Multi-company ridesharing}\label{multi}
\revtwo{
In the same New York City setting, we consider the case of multiple ridesharing companies mentioned in Section~\ref{distri}, e.g. 2, sharing the whole customer set in Manhattan, but with a twist. We consider the same situation of users booking via the proprietary applications of the companies, and a central agent, e.g. the city authority, ensuring optimal assignment of the requests. The more cars the companies have on the roads, the more they are assigned to the users, and the more market share they will get (as they guarantee better quality of service, and gain reputation). Hence, there is a negative competition taking place bringing more and more cars on the roads, which in turn augments traffic congestion and pollution and is detrimental to the well-being of the citizens. Therefore, here we investigate the situation in which the central agent allocates the customers to the companies based not only on a global welfare cost (our cost function TD) \textit{but} also on a predefined market share, for instance proportional to the number of issued licenses for the companies. Note that the discussion presented in this section is particularly timely considering the recent context: New York city council became in August 2018 the first US city to halt new licenses for ridesharing services, see e.g.~\cite{uber}, and cities are looking at mechanisms to cap the number of issued licenses in a reasonable fashion.}

In this scenario, the two companies have defined market share (say $25\%$ and $75\%$) and the city authority enforces it via carefully engineering the cost function, so that the two companies expect to receive $25\%$ and $75\%$ of the city customers, respectively. 

The cost function is modified adding a penalty term to each $c_{ij}$ terms as
\begin{eqnarray}
c_{ij}' =& c_{ij} + (s_1 - q (s_1 + s_2))^3, & \textrm{if } i \in \textrm{ company 1}, \\
c_{ij}' =& c_{ij} + (s_2 - (1-q) (s_1 + s_2))^3, & \textrm{if } i \in \textrm{ company 2}, 
\end{eqnarray}
where $s_1$ and $s_2$ are the cumulative served customers from $t=0$ till the current time, for company 1 and company 2, respectively, while $q\in[0,1]$ is the market share of company 1. The added penalty serves as a integral feedback loop that increases the appeal of giving customers to company $i$ if its shares are lower than agreed, and vice-versa. We further define violation as $V = s_2 - (1-q) (s_1 + s_2)$. 

\revtwo{
We consider company $1$ having $100$ vehicles with $75\%$ of market share, while company $2$ having $50$ vehicles and with $25\%$ of market share. Naturally, if the vehicles of each of the two companies were identically spatially-distributed and if they were designed to have $66.7\%$ and $33.3\%$ of the market shares, respectively, we expect the standard algorithm with cost function TD to have a quality of service comparable to the centralised case. Here we explore the case in which the companies size their fleet differently (to the number of issued licenses) and are not aware of the centralized best. 
}

\begin{figure}
\centering
\includegraphics[width=0.6\textwidth]{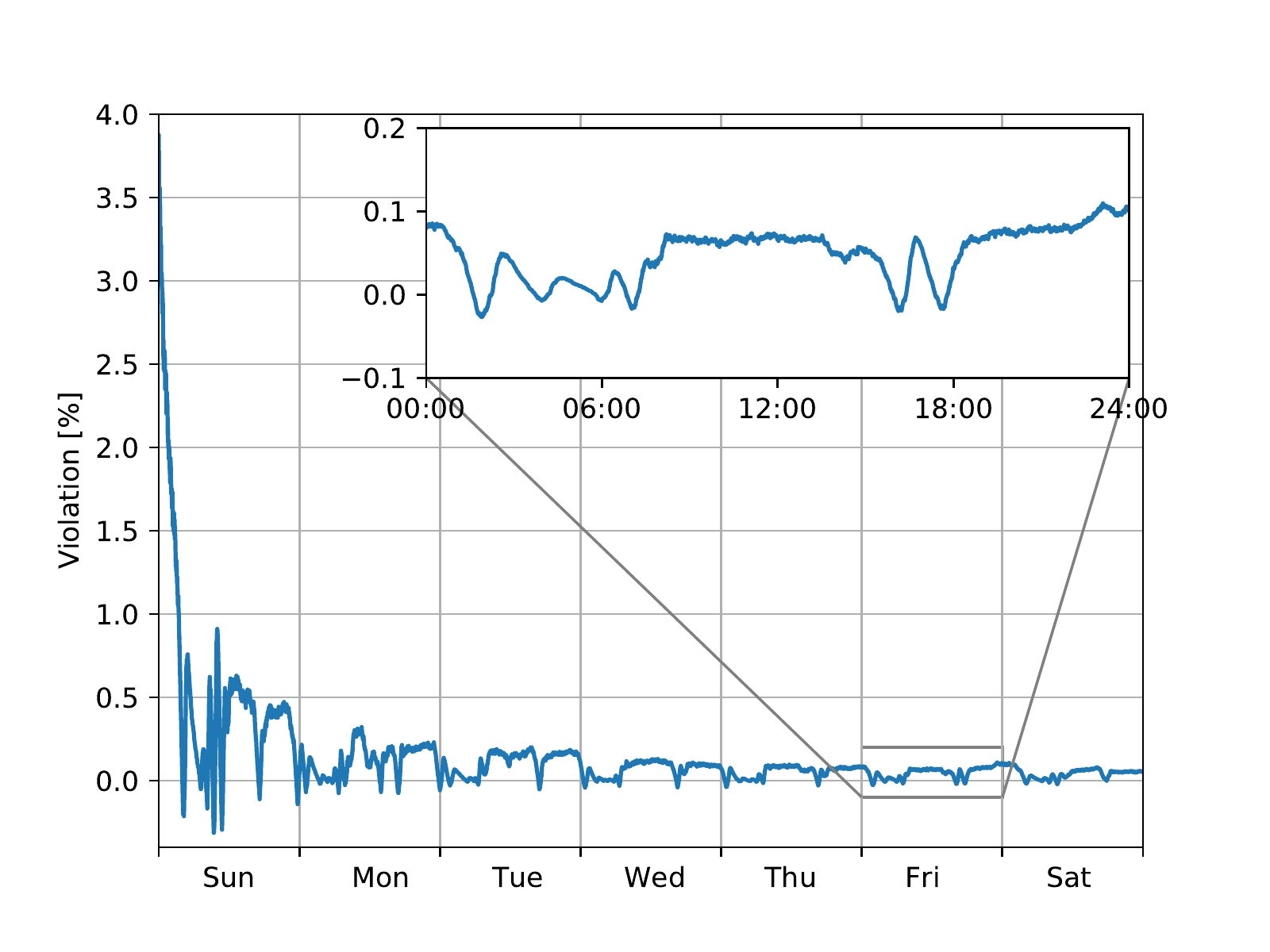}
\caption{Violation of the defined market shares. A positive violation means that company 2 gets more customers than agreed and vice-versa. }
\label{fig:violation}
\end{figure} 

Figure~\ref{fig:violation} reports how the violation $V$ (expressed in percentage) evolves over time. A positive violation means that company 2 receives more cars than agreed, and vice-versa. As one can see, the violation settles around $\pm 0.1\%$ from the second day of service, indicating that the modified cost TD' does the job.

\begin{table}[H]
\footnotesize
\centering
\caption{Comparison between optimal and market share solutions;  market shares: $75\%$ and $25\%$, respectively.}
\label{tab:ny-r-c}
\begin{tabular}{ccccccccccccc}
\toprule
vehicles & vehicles & customers & $c$ & max$n$ & method & $h$ & SR & waiting & waiting  &  detour  & detour  \\ 
 & [\%] & [\%] & &  & &[s] & [\%] & y [min] & n [min]  &  y [min] & n [min]   \\ 
\toprule
\rowcolor{Gray}
150 & 5.0 & 5.0 & 4 & 8 & TD & 10 & 95.75 & 3.54 & 3.43 & 2.53 & 2.53 \\ 
100+50 & 5.0 & 5.0 & 4 & 16 & TD' & 10 & 89.49 & 4.83 & 4.72 & 4.76 & 4.76\\ 
\toprule
\end{tabular}
\end{table} 

Table~\ref{tab:ny-r-c} reports the performance of the scenario with respect to a monopolistic setting. As we see, performance deteriorates on all the metrics (even if we augment max$n$ to 16). In order to obtain a similar level of performance to the monopolistic case, company 1 will have to increase its fleet size, because company 2 will most likely not reduce its fleet, unless enforced to. As a result, the total fleet size will increase, as well as the number of vehicles on the roads and traffic congestion. 

\revtwo{An optimal assignment strategy will tend to reward the companies with the most vehicles on the roads. This is not desired as companies will aim at increasing their active fleet size to increase their market share, hence bringing more congestion and pollution on the roads. Although an intuitive idea, the engineering of the cost function to ensure a fair distribution of the users across the companies is not a good one. Indeed we have showed that it can decrease the quality of the service provided in a significant way. Instead, the city authorities must look into adapting the fleet size of the companies to the time-varying demand (number of requests per time window). This underlines the needs for strong policies by city authorities to regulate the number of cars utilized by ridesharing companies. Further work on the topic falls beyond the scope of this paper, but we believe is of great interest.}

\rev{
\subsection{Melbourne Metropolitan Area simulations}

In this section, we evaluate the proposed ridesharing solution against a different dataset, that is available at the repository \href{url}{https://github.com/davidrey123/Ridesharing}. This dataset contains instances of realistic ridesharing demand scenarios corresponding to 0.25\%, 0.5\% and 0.75\% of the total demand for the Melbourne Metropolitan area, Australia, on a standard day between 06:00 hours and 21:00 hours. \revtwo{The percentage may seem small with respect to the total demand but one should note that only the demand related to ridesharing is considered here, as many travellers do not fall in this category. The interested reader can refer to~\citep{davidtorche} for more detail on the dataset. This dataset was indeed already utilized in a previous work for the peer-to-peer ridesharing problem~\citep{davidtorche}.} 

We focus on a representative ridesharing scenario, i.e., 0.5\%, instance $M1$, and the same sampling period as in~\cite{davidtorche}, i.e., 2 min. It is interesting to note that the requests are now \textit{scheduled} requests, which means that the earliest and latest pick-up times are defined by the customers themselves, differently to the New York taxi dataset which focus on \textit{instantaneous} requests, where a maximum waiting time (or delay) of 7 min was defined. This enables to test the behavior of our proposed ridesharing solution in a different framework, as well as in a different geographical area. 

As discussed in Section~\ref{sec:vehicle}, we have assumed so far that a vehicle with capacity $C$ can have at most 2$C$ requests in its pipeline of scheduled pickup/delivery events. Given the larger time horizon we are considering with the Melbourne dataset, we enlarge this to consider at most 4$C$ requests. In the rest of the section, we present results for $C=4$, and a maximum number of vehicles considered in the context mapping module of max$n=10$, see Section~\ref{sec:context-mapping}.

\begin{figure}
\centering
\includegraphics[width=0.50\textwidth]{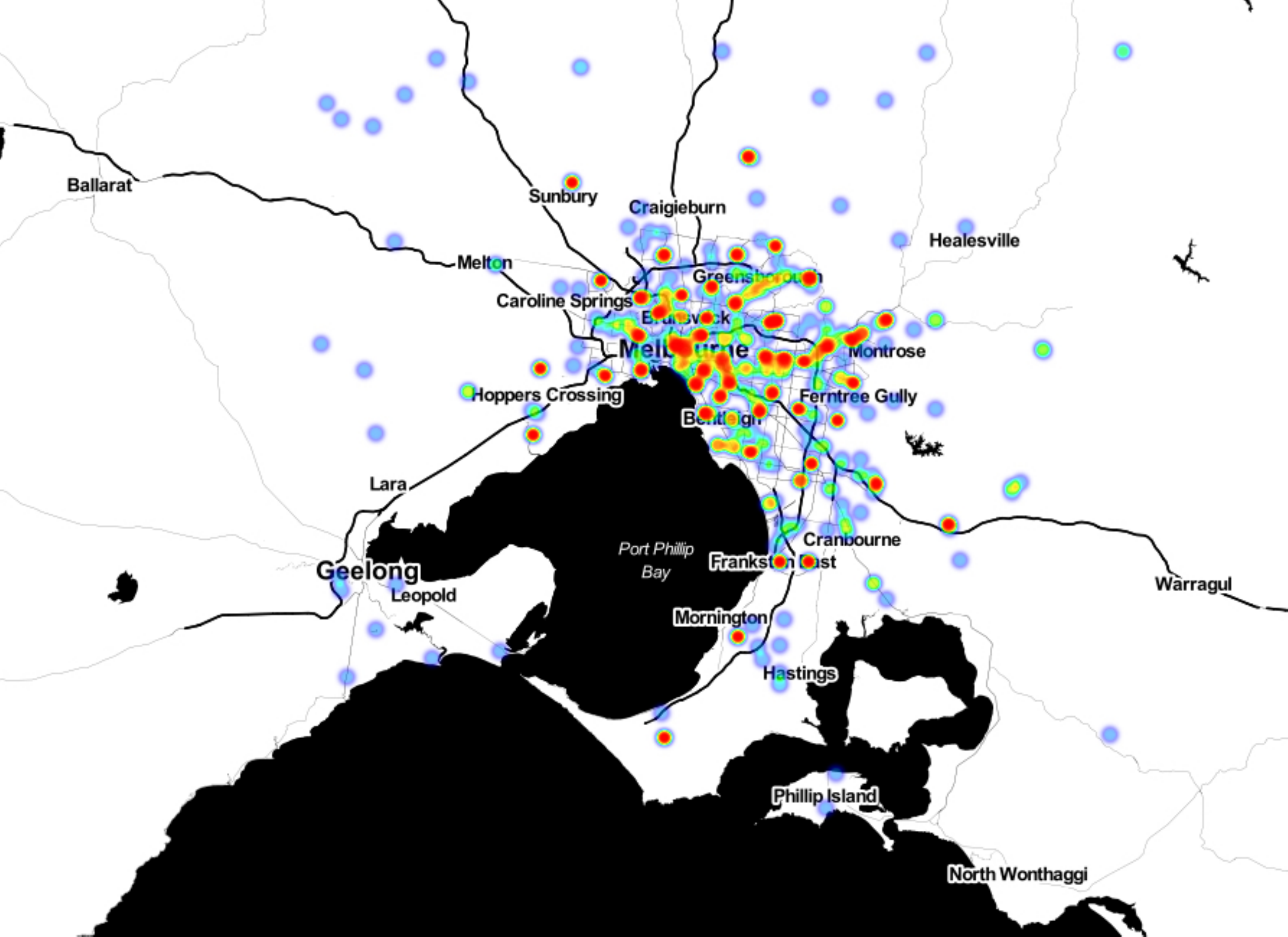}
\caption{A snapshot of the Melbourne Metropolitan Area with a fleet of 1000 vehicles with capacity of 4 at 21:00 hours; the heat map represents how many customers are in the vehicles at that moment.}
\label{fig.sim-1M}
\end{figure}

Figure~\ref{fig.sim-1M} illustrates the density of customers in the vehicles at the end of the simulation experiment. In analogy with the results for NYC, it is interesting to note that the pool of requests is more spread out in space, which intuitively limits the usefulness of ridesharing applications. We shall see that ridesharing is still attractive despite this configuration.

Figures~\ref{fig.sim-2M} and~\ref{fig.sim-3M} present the evolution of the detour times per customer and of the average number of passengers per vehicle in the time frame under consideration for $600$, $800$, and $1000$ vehicles. 
The morning and evening peak of demand are clearly visible in the two figures, particularly for the case of $800$ and $1000$ vehicles. Indeed, Figure~\ref{fig.sim-3M} shows that less than 1 vehicle per passenger is enough in the non-peak hours, i.e. ridesharing may not be needed, but this changes at peak time. On the contrary, the situation with $600$ vehicles shows a permanent use of the ridesharing offering, on average. 

Regarding the computational complexity of the algorithm, near-real time is achieved as indicated in Figure~\ref{fig.sim-4M}. Note that this result applies despite the fact that, contrary to the New York dataset, the travel times $\tau$ cannot be precomputed but have to be computed on-the-fly with OSRM due to the large network size. Note that the average computational time per $2$ min period is even of $1.9$ min in the $600$ vehicles case, as presented in Table~\ref{tab.2}. Regarding the service rate, it already reaches $75.7\%$ in the $600$ vehicles scenario to attain $100\%$ for $1000$ vehicles.

\begin{figure}
\centering
\includegraphics[width=0.99\textwidth]{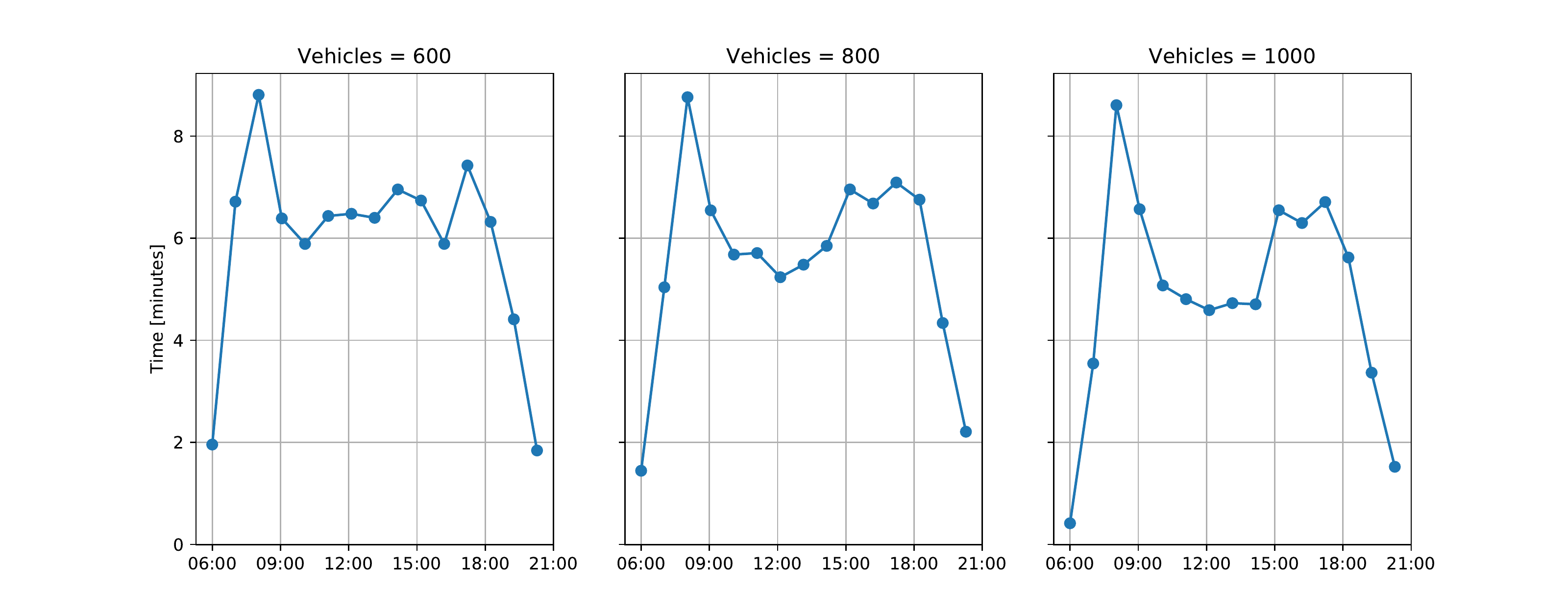}
\caption{Detour times per customer in the considered time frame.}
\label{fig.sim-2M}
\end{figure}

\begin{figure}
\centering
\includegraphics[width=0.99\textwidth]{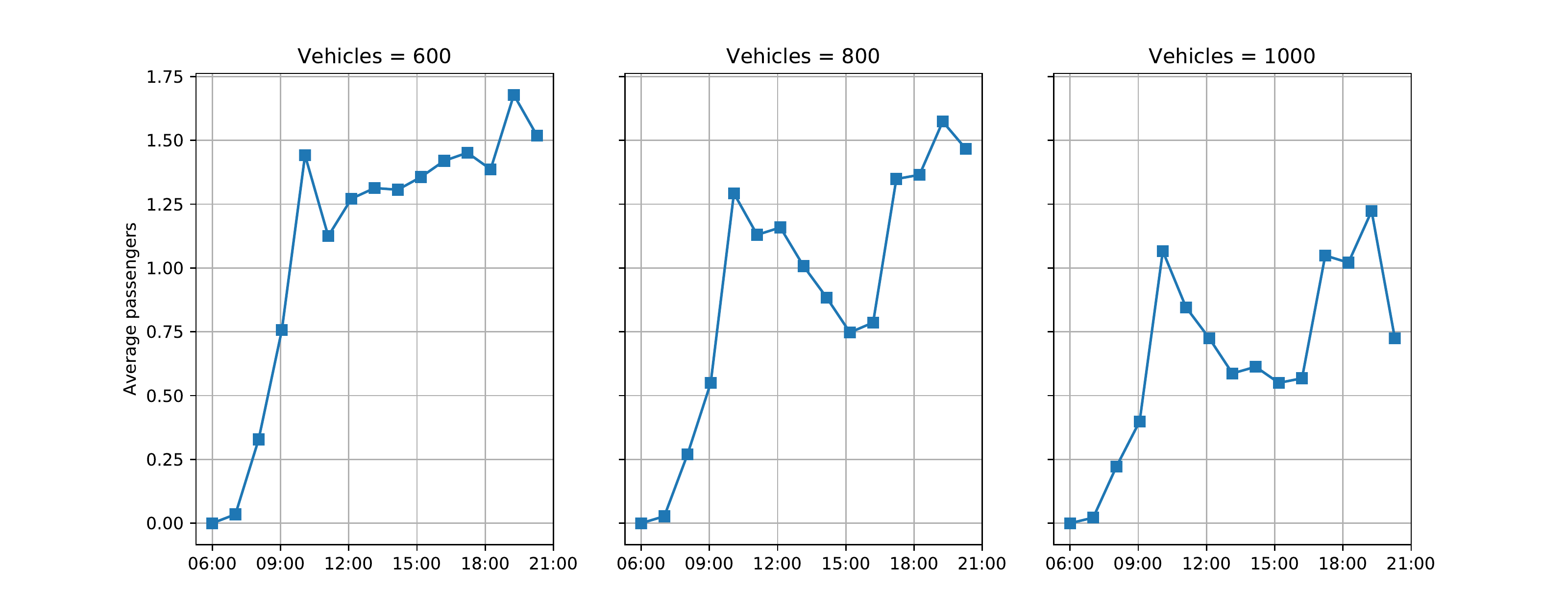}
\caption{Average number of passengers per vehicle in the considered time frame.}
\label{fig.sim-3M}
\end{figure}

\begin{table}
	\footnotesize
	\centering
	\caption{Service rate, detour and computational times for $600$, $800$ and $1000$ vehicles.}
	\label{tab.2}
	\begin{tabular}{cccccccccc}
		\toprule
		vehicles & customers & $c$ & max$n$ & cost & $h$ & SR &   detour  & detour & comp. time  \\ 
		&  [\%] & &  & &[min] & [\%] &   y [min] & n [min] & [min]  \\ 
		\toprule
		600 &  100.0 & 4 & 10 & TD & 2 & 75.68 &  5.90 & 5.91 & 1.99 \\ 
		\rowcolor{Gray}
		800 & 100.0 & 4 & 10 & TD & 2 & 96.06 & 5.58 & 5.58  & 2.61 \\ 
		1000 & 100.0 & 4 & 10 & TD & 2 & 100.00 & 4.87 & 4.87 & 2.96  \\  
		\toprule
	\end{tabular}
\end{table}

\begin{figure}
\centering
\includegraphics[width=0.99\textwidth]{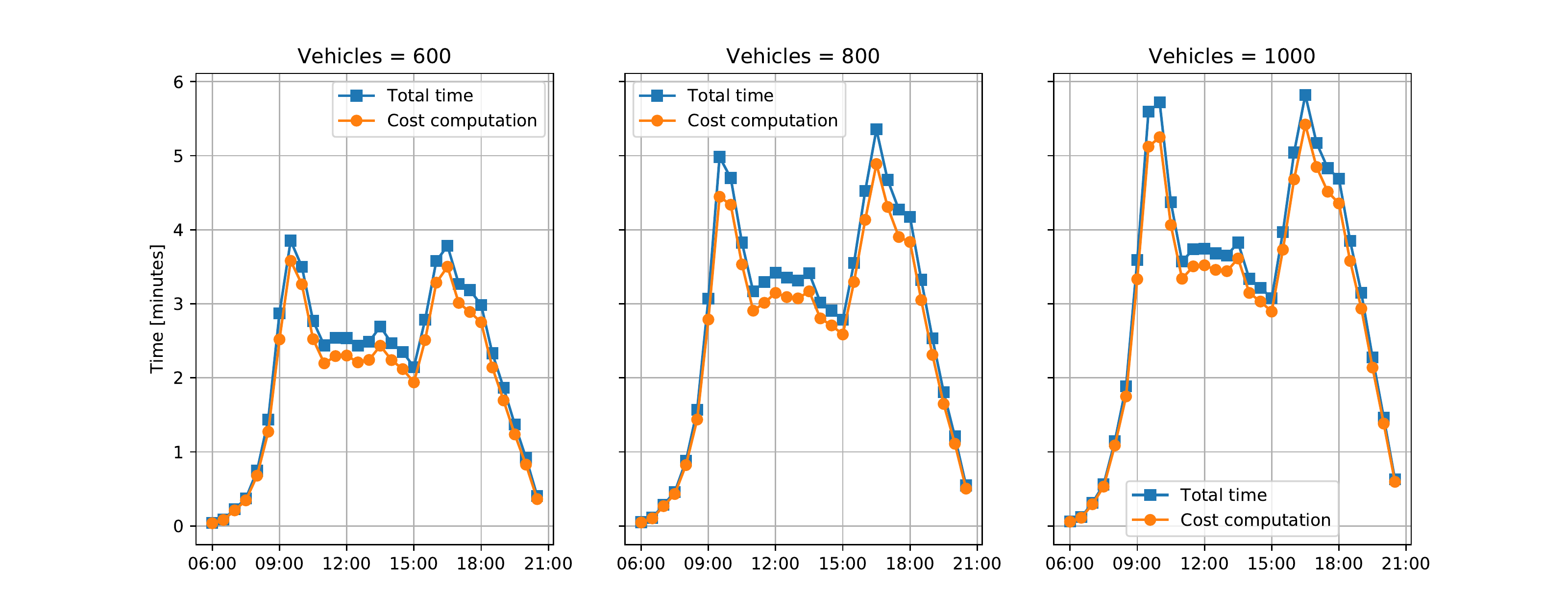}
\caption{Computational time per call of the ridesharing service in the considered time frame.}
\label{fig.sim-4M}
\end{figure} 

Finally, in the light of the discussion in the previous work~\citep{davidtorche}, which focused on the peer-to-peer ridesharing problem, i.e., the matching between riders and drivers, we can make the following observations. In~\cite{davidtorche}, for instance $M1$, on the global run, a total of $25550$ drivers are attempted to be matched with a total of $20250$ riders, that is on average $57$ vehicles available and $45$ requests per batch of $2$ min. As the average trip duration is observed to be of $16$ min, therefore $8$ batch periods, it would take respectively $600/8=75$ and $800/8=100$ vehicles per batch of $2$ min to satisfy the demand (up to respectively $76\%$ and $96\%$) on average considering on-demand ridesharing, which is an explanation of the low matching rate presented in~\cite{davidtorche}. Besides, there are $8700$ taxis available in the Metropolitan Melbourne area~\citep{Melbournetaxi}, so it interesting to note that less than approximately $10\%$ of them would be needed to meet $96\%$ of the requests of the considered instance. These preliminary results certainly require further investigation, however we believe that similar city-wise numerical experiments should be conducted as they shed some light on limiting the number of vehicles on the roads to meet the demand.   

} \color{black}


\section{Conclusions and open research questions}\label{sec:open}

We have presented a novel, computational efficient, dynamic ridesharing algorithm, based on a federated optimization architecture and a suitably defined linear assignment problem. With the help of numerical simulations based on the New York Taxi dataset \rev{and on a dataset coming from Melbourne Metropolitan Area}, we have shown the performance of the algorithm and its comparison with state-of-the-art approaches. This paper further advocates that small and medium enterprises with small market shares ($5$-$10\%$) can achieve similar benefits to the ones offered by a more global ridesharing service. Naturally, such benefits can only be obtained if the number of vehicles in circulation at the scale of a city are \rev{carefully planned, as we have argued via simulations featuring two ridesharing companies that buy their market share upfront}.  

Several open questions are left to future research works. First, the use of demand forecasting in scheduling and rebalancing could significantly improve quality of service, particularly in low market share regimes ($< 5\%$). Second, the interplay between intelligent public transportation, e.g., ridesharing buses with high capacity, dedicated ridesharing systems, and traditional transport could shape new mobility paradigms. In this context, multi-legged and multi-modal ridesharing are of particular interest. \rev{Finally, studying competition among ridesharing companies as well as incentives and penalties is crucial to make the best of the emerging mobility-on-demand world.}

\newpage
\footnotesize
\bibliographystyle{agsm}
\bibliography{Autopilot_Literature}

\end{document}